\documentclass[11pt, reqno]{article}
\usepackage{amsfonts}
\usepackage{html}
\usepackage{url}
\usepackage[dcucite]{harvard}
\usepackage{graphicx,amsmath,amssymb,epsfig,harvard}
\usepackage{colortbl}

\renewcommand{\thepage}{}

\input{tcilatex}
\begin{document}

\title{Local Identification of Nonparametric and Semiparametric Models%
\thanks{%
The National Science Foundation (Grants SES-0838161 and SES-1132399),
European Research Council (ERC-2009-StG-240910-ROMETA), and the National
Research Foundation of Korea (NRF-2011-327-B00073) provided financial
support for this paper. Helpful comments were provided by a co-editor, three
anonymous referees, D. Andrews, D. Chetverikov, K. Evdokimov, J.P. Florens,
L. Hansen, J. Heckman, S. Hoderlein, H. Ichimura, T. Komarova, O. Linton, A.
Santos and participants in seminars at June 2010 Cowles Conference, June
2010 Beijing Conference, June 2010 LSE Conference, and December 2010 EC2
conference.}}
\author{Xiaohong Chen \\
Department of Economics \\
Yale \and Victor Chernozhukov \\
Department of Economics\\
MIT \and Sokbae Lee \\
Department of Economics\\
Seoul National University \and Whitney K. Newey \\
Department of Economics\\
MIT}
\date{ 2009\\
Revised, April 2013}
\maketitle

In parametric, nonlinear structural models a classical sufficient condition
for local identification, like Fisher (1966) and Rothenberg (1971), is that
the vector of moment conditions is differentiable at the true parameter with
full rank derivative matrix. We derive an analogous result for the \textit{%
nonparametric}, nonlinear structural models, establishing conditions under
which an infinite-dimensional analog of the full rank condition is
sufficient for local identification. Importantly, we show that additional
conditions are often needed in nonlinear, nonparametric models to avoid
nonlinearities overwhelming linear effects. We give restrictions on a
neighborhood of the true value that are sufficient for local identification.
We apply these results to obtain new, primitive identification conditions in
several important models, including nonseparable quantile instrumental
variable (IV) models, single-index IV models, and semiparametric
consumption-based asset pricing models. \newline

\textbf{JEL Classification: }C12, C13, C23

\textbf{Keywords:}\textit{\ }Identification, Local Identification,
Nonparametric Models, Asset Pricing.

\baselineskip=22pt\newpage \setcounter{page}{1}\renewcommand{\thepage}{[%
\arabic{page}]}\renewcommand{\theequation}{\arabic{section}.%
\arabic{equation}}

\section{Introduction}

There are many important models in econometrics that give rise to
conditional moment restrictions. These restrictions often take the form%
\begin{equation*}
\text{E}[\rho (Y,X,\alpha _{0})|W]=0\text{,}
\end{equation*}%
where $\rho (Y,X,\alpha )$ has a known functional form but $\alpha _{0}$ is
unknown. Parametric models (i.e., models when $\alpha _{0}$\ is finite
dimensional) of this form are well known from the work of Hansen (1982),
Chamberlain (1987), and others. Nonparametric versions (i.e., models when $%
\alpha _{0}$\ is infinite dimensional) are motivated by the desire to relax
functional form restrictions. Identification and estimation of linear
nonparametric conditional moment models have been studied by Newey and
Powell (2003), Hall and Horowitz (2005), Blundell, Chen, and Kristensen
(2007), Darolles, Fan, Florens, and Renault (2011), and others.

The purpose of this paper is to derive identification conditions for $\alpha
_{0}$ when $\rho $ may be nonlinear in $\alpha $ and for other nonlinear
nonparametric models. Nonlinear models are important. They include models
with conditional quantile restrictions, as discussed in Chernozhukov and
Hansen (2005) and Chernozhukov, Imbens, and Newey (2007), and various
economic structural and semiparametric models, as further discussed below.
In this paper we focus on conditions for local identification of these
models. It may be possible to extend these results to provide global
identification conditions.

In parametric models there are easily interpretable rank conditions for
local identification, as shown in Fisher (1966) and Rothenberg (1971). We
give a pair of conditions that are sufficient for parametric local
identification from solving a set of equations. They are a) pointwise
differentiability at the true value, and b) the rank of the derivative
matrix is equal to the dimension of the parameter $\alpha _{0}$. We find
that the nonparametric case is different. Differentiability and the
nonparametric version of the rank condition may not be sufficient for local
identification. We suggest a restriction on the neighborhood that does give
local identification, via a link between curvature and an identification
set. We also give more primitive conditions for Hilbert spaces, that include
interesting econometric examples. In addition we consider semiparametric
models, providing conditions for identification of a finite dimensional
Euclidean parameter. These conditions are based on "partialling out" the
nonparametric part and allow for identification of the parametric part even
when the nonparametric part is not identified.

The usefulness of these conditions is illustrated by three examples. One
example gives primitive conditions for local identification of the
nonparametric endogenous quantile models, where primitive identification
conditions had only been given previously for discrete regressors. Another
example gives conditions for local identification of a semiparametric index
model with endogeneity. There we give conditions for identification of
parametric components when nonparametric components are not identified. The
third example gives sufficient conditions for local identification of a
semiparametric consumption capital asset pricing model.

In relation to previous literature, in some cases the nonparametric rank
condition is a local version of identification conditions for linear
conditional moment restriction models that were considered in Newey and
Powell (2003). Chernozhukov, Imbens, and Newey (2007) also suggested
differentiability and a rank condition for local identification but did not
recognize the need for additional restrictions on the neighborhood. Florens
and Sbai (2010) gave local identification conditions for games but their
conditions do not apply to the kind of conditional moment restrictions that
arise in instrumental variable settings and are a primary subject of this
paper.

Section 2 presents general nonparametric local identification results and
relates them to sufficient conditions for identification in parametric
models. Section 3 gives more primitive conditions for Hilbert spaces and
applies them to the nonparametric endogenous quantile model. Section 4
provides conditions for identification in semiparametric models and applies
these to the endogenous index model. Section 5 discusses the semiparametric
asset pricing example and Section 6 briefly concludes. The Appendix contains
additional lemmas and all of the proofs.

\section{Nonparametric Models}

\subsection{The Setting and Definition of Local Identification}

To help explain the nonparametric results and give them context we give a
brief description of sufficient conditions for local identification in
parametric models. Let $\alpha $ be a $p\times 1$ vector of parameters and $%
m(\alpha )$ a $J\times 1$ vector of functions with $m(\alpha _{0})=0$ for
the true value $\alpha _{0}$. Also let $\left\vert \cdot \right\vert $
denote the Euclidean norm in either $\mathbb{R} ^{p}$ or $\mathbb{R} ^{J}$
depending on the context. We say that $\alpha _{0}$ is locally identified if
there is a neighborhood of $\alpha _{0}$ such that $m(\alpha )\neq 0$ for
all $\alpha \neq \alpha _{0}$ in the neighborhood. Let $m^{\prime }$ denote
the derivative of $m(\alpha )$ at $\alpha _{0}$ when it exists. Sufficient
conditions for local identification can be stated as follows:

\bigskip

\textit{If }$m(\alpha )$\textit{\ is differentiable at }$\alpha _{0}$\textit{%
\ and }$rank(m^{\prime })=p$\textit{\ then }$\alpha _{0}$\textit{\ is
locally identified.}

\bigskip

This statement is proved in the Appendix. Here the sufficient conditions for
parametric local identification are pointwise differentiability at the true
value $\alpha _{0}$ and the rank of the derivative equal to the number of
parameters.

In order to extend these conditions to the nonparametric case we need to
modify the notation and introduce structure for infinite dimensional spaces.
Let $\alpha $ denote a function with true value $\alpha _{0}$ and $m(\alpha
) $ a function of $\alpha $, with $m(\alpha _{0})=0$. Conditional moment
restrictions are an important example where $\rho (Y,X,\alpha )$ is a finite
dimensional residual vector depending on an unknown function $\alpha $ and $%
m(\alpha )=\text{E}[\rho (Y,X,\alpha )|W]$. We impose some mathematical
structure by assuming that $\alpha \in \mathcal{A}$, a Banach space with
norm $\left\Vert \cdot \right\Vert _{\mathcal{A}}$ and $m(\alpha )\in 
\mathcal{B}$, a Banach space with a norm $\left\Vert \cdot \right\Vert _{%
\mathcal{B}}$, i.e. $m:\mathcal{A}\mapsto \mathcal{B}$. The restriction of
the model is that $\left\Vert m(\alpha _{0})\right\Vert _{\mathcal{B}}=0$.
The notion of local identification we consider is:

\bigskip

\textsc{Definition:} $\alpha _{0}$\textit{\ is locally identified on }$%
\mathcal{N}\subseteq \mathcal{A}$\textit{\ if }$\left\Vert m(\alpha
)\right\Vert _{\mathcal{B}}>0$\textit{\ for all }$\alpha \in \mathcal{N}$ 
\textit{with }$\alpha \neq \alpha _{0}.$

\bigskip

This local identification concept is more general than the one introduced by
Chernozhukov, Imbens and Newey (2007). Note that local identification is
defined on a set $\mathcal{N}$ in $\mathcal{A}$. Often there exists an $%
\varepsilon >0$ such that $\mathcal{N}$ is a subset of an open ball 
\begin{equation*}
\mathcal{N}_{\varepsilon }\equiv \{\alpha \in \mathcal{A}:\left\Vert \alpha
-\alpha _{0}\right\Vert _{\mathcal{A}}<\varepsilon \}.
\end{equation*}%
It turns out that it may be necessary for $\mathcal{N}$ to be strictly
smaller than an open ball $\mathcal{N}_{\varepsilon }$ in $\mathcal{A}$, as
discussed below.

\subsection{Local Identification via Full-Rank Conditions}

The nonparametric version of the derivative will be a bounded (i.e.,
continuous) linear map $m^{\prime }:\mathcal{A}\mapsto \mathcal{B}$. Under
the conditions we give, $m^{\prime }$ will be a G\^{a}teaux derivative at $%
\alpha _{0},$ that can be calculated as%
\begin{equation}
m^{\prime }h=\frac{\partial }{\partial t}m(\alpha _{0}+th)|_{t=0}
\label{Gauteaux}
\end{equation}%
for $h\in \mathcal{A}$ and $t$ a scalar. Sometimes we also require that for
any $\delta >0$ there is $\varepsilon >0$ with 
\begin{equation*}
\frac{\left\Vert m(\alpha )-m(\alpha _{0})-m^{\prime }(\alpha -\alpha
_{0})\right\Vert _{\mathcal{B}}}{\left\Vert \alpha -\alpha _{0}\right\Vert _{%
\mathcal{A}}}<\delta \text{.}
\end{equation*}%
for all $\alpha \in \mathcal{N}_{\varepsilon }$. This is Fr\'{e}chet
differentiability of $m(\alpha )$ at $\alpha _{0}$ (which implies that the
linear map $m^{\prime }:\mathcal{A}\mapsto \mathcal{B}$ is continuous). Fr%
\'{e}chet differentiability of estimators that are functionals of the
empirical distribution is known to be too strong, but\ is typically
satisfied in local identification analysis, as shown by our examples.

In parametric models the rank condition is equivalent to the null space of
the derivative matrix being zero. The analogous nonparametric condition is
that the null space of the linear map $m^{\prime }$ is zero, as follows:

\bigskip

\textsc{Assumption 1 (Rank Condition):} \textit{There is a set $\mathcal{N}%
^{\prime }$ such that }$\left\Vert m^{\prime }(\alpha -\alpha
_{0})\right\Vert _{\mathcal{B}}>0$\textit{\ for all }$\alpha \in $\textit{$%
\mathcal{N}^{\prime }$ with }$\alpha \neq \alpha _{0}.$

\bigskip

This condition is familiar from identification of a linear conditional
moment model where $Y=\alpha _{0}(X)+U$ and $\text{E}[U|W]=0.$ Here $\rho
(Y,X,\alpha )=Y-\alpha (X)$, so that $m(\alpha )=\text{E}[Y-\alpha (X)|W]$
and $m^{\prime }h=-\text{E}[h(X)|W]$. In this case Assumption 1 requires
that $\text{E}[\alpha (X)-\alpha _{0}(X)|W]\neq 0$ for any $\alpha \in $%
\textit{$\mathcal{N}^{\prime }$ }with $\alpha -\alpha _{0}\neq 0$. For $%
\mathcal{N}^{\prime }=\mathcal{A}$ this is the completeness condition
discussed in Newey and Powell (2003). Andrews (2011) has recently shown that
if $X$ and $W$ are continuously distributed, there are at least as many
instruments in $W$ as regressors in $X$, and the conditional distribution of 
$X$ given $W$ is unrestricted (except for a mild regularity condition), then
the completeness condition holds generically, in a sense defined in that
paper. In Section 3 we also give a genericity result for a different range
of models. For this reason we think of Assumption 1 with $\mathcal{N}%
^{\prime }=\mathcal{A}$ as a weak condition when there are as many
continuous instruments\textbf{\ }$W$\ as the endogenous regressors $X$, just
as it is in a parametric linear instrumental variables model with
unrestricted reduced form. It is also an even weaker condition if some
conditions are imposed on the deviations, so in the statement of Assumption
1 we allow it to hold only on $\mathcal{N}^{\prime }\subset \mathcal{A}.$
For example, if we restrict $\alpha -\alpha _{0}$ to be a bounded function
of $X$, then in linear conditional moment restriction models Assumption 1
only requires that the conditional distribution of $X$ given $W$ be bounded
complete, which is known to hold for even more distributions than does
completeness. This makes Assumption 1 even more plausible in models where $%
\alpha _{0}$ is restricted to be bounded, such as\textbf{\ }in Blundell,
Chen and Kristensen (2007). See, for example, Mattner (1993), Chernozhukov
and Hansen (2005), D'Haultfoeuille (2011), and Andrews (2011) for
discussions of completeness and bounded completeness.

Fr\'{e}chet differentiability and the rank condition are not sufficient for
local identification in an open ball $\mathcal{N}_{\varepsilon }$ around $%
\alpha _{0},$ as we further explain below. One condition that can be added
to obtain local identification is that $m^{\prime }:\mathcal{A}\mapsto 
\mathcal{B}$ is onto.

\bigskip

\textsc{Theorem 1: }\textit{If} $m(\alpha )$\textit{\ is Fr\'{e}chet
differentiable at }$\alpha _{0},$\textit{\ the rank condition is satisfied
on $\mathcal{N}^{\prime }$ }$=\mathcal{N}_{\varepsilon }$ \textit{for some }$%
\varepsilon >0$, \textit{and }$m^{\prime }:\mathcal{A}\mapsto \mathcal{B}$ 
\textit{is onto, then }$\alpha _{0}$\textit{\ is locally identified on }$%
\mathcal{N}_{\tilde{\varepsilon}}$ \textit{for some} $\tilde{\varepsilon}$ 
\textit{with }$0<\tilde{\varepsilon}\leq \varepsilon .$

\bigskip

This result extends previous nonparametric local identification results by
only requiring pointwise Fr\'{e}chet differentiability at $\alpha _{0},$
rather than continuous Fr\'{e}chet differentiability in a neighborhood of $%
\alpha _{0}$. This extension may be helpful for showing local identification
in nonparametric models, because conditions for pointwise Fr\'{e}chet
differentiability are simpler than for continuous differentiability in
nonparametric models.

Unfortunately, the assumption that $m^{\prime }$ is onto is too strong for
many econometric models, including many nonparametric conditional moment
restrictions. An onto $m^{\prime }$ implies that $m^{\prime }$ has a
continuous inverse, by the Banach Inverse Theorem (Luenberger, 1969, p.
149). The inverse of $m^{\prime }$ may not be continuous for nonparametric
conditional moment restrictions, as discussed in Newey and Powell (2003).
Indeed, the discontinuity of the inverse of $m^{\prime }$ is a now well
known ill-posed inverse problem that has received much attention in the
econometrics literature, e.g. see the survey of Carrasco, Florens, and
Renault (2007). Thus, in many important econometric models Theorem 1 cannot
be applied to obtain local identification.

It turns out that $\alpha _{0}$ may not be locally identified on any open
ball in ill-posed inverse problems, as we show in an example below. The
problem is that for infinite dimensional spaces $m^{\prime }(\alpha -\alpha
_{0})$ may be small when $\alpha -\alpha _{0}$ is large. Consequently, the
effect of nonlinearity, that is related to the size of $\alpha -\alpha _{0}$%
, may overwhelm the identifying effect of nonzero $m^{\prime }(\alpha
-\alpha _{0})$, resulting in $m(\alpha )$ being zero for $\alpha $ close to $%
\alpha _{0}$.

We approach this problem by restricting the deviations $\alpha -\alpha _{0}$
to be small when $m^{\prime }(\alpha -\alpha _{0})$ is small. The
restrictions on the deviations will be related to the nonlinearity of $%
m(\alpha )$ via the following condition:

\bigskip

\textsc{Assumption 2: }\textit{There are }$L\geq 0,r\geq 1$\textbf{\ }%
\textit{and a set $\mathcal{N}^{\prime \prime }$ such that for all }$\alpha
\in $\textit{\ $\mathcal{N}^{\prime \prime },$} 
\begin{equation*}
\left\Vert m(\alpha )-m(\alpha _{0})-m^{\prime }(\alpha -\alpha
_{0})\right\Vert _{\mathcal{B}}\leq L\left\Vert \alpha -\alpha
_{0}\right\Vert _{\mathcal{A}}^{r}.
\end{equation*}

\bigskip

This condition is general. It includes the linear case where \textit{$%
\mathcal{N}^{\prime \prime }=\mathcal{A}$ }and $L=0$. It also includes Fr%
\'{e}chet differentiability, where $r=1,$ $L$ is any positive number and $%
\mathcal{N}^{\prime \prime }=\mathcal{N}_{\varepsilon }$\textit{\ }for any
sufficiently small $\varepsilon >0$. Cases with $r>1,$ are analogous to H%
\"{o}lder continuity of the derivative in finite dimensional spaces, with $%
r=2$ corresponding to twice continuous Fr\'{e}chet differentiability. We
would only have $r>2$ when the second derivative is zero. This condition
applies to many interesting examples, as we will show in the rest of the
paper. The term $L\left\Vert \alpha -\alpha _{0}\right\Vert _{\mathcal{A}%
}^{r}$ represents a magnitude of nonlinearity that is allowed for $\alpha
\in $\textit{$\mathcal{N}^{\prime \prime }$. }The following result uses
Assumption 2 to specify restrictions on $\alpha $ that are sufficient for
local identification.

\bigskip

\textsc{Theorem 2: }\textit{If Assumption 2 is satisfied then }$\alpha _{0}$ 
\textit{is locally identified on $\mathcal{N=N}^{\prime \prime }\cap 
\mathcal{N}^{\prime \prime \prime }$ with $\mathcal{N}^{\prime \prime \prime
}=$}$\{\alpha :\left\Vert m^{\prime }(\alpha -\alpha _{0})\right\Vert _{%
\mathcal{B}}>L\left\Vert \alpha -\alpha _{0}\right\Vert _{\mathcal{A}}^{r}\}$%
.

\bigskip

The strict inequality in \textit{$\mathcal{N}^{\prime \prime \prime }$ }is
important for the result. It does exclude $\alpha _{0}$ from \textit{$%
\mathcal{N}$}, but that works because local identification specifies what
happens when $\alpha \neq \alpha _{0}$. This result includes the linear
case, where $L=0,$ $\mathit{\mathcal{N}}^{\prime \prime }=\mathcal{A}$, and $%
\mathcal{N}=\mathcal{N^{\prime \prime \prime }}=\mathcal{N}^{\prime }$. It
also includes nonlinear cases where only Fr\'{e}chet differentiability is
imposed, with $r=1$ and $L$ equal to any positive constant. In that case $%
\mathcal{N}^{\prime \prime }=\mathcal{N}_{\varepsilon }$ for some $%
\varepsilon $ small enough and $\alpha \in \mathcal{N}^{\prime \prime \prime
}$ restricts $\alpha -\alpha _{0}$ to a set where the inverse of $m^{\prime
} $ is continuous by requiring that $\left\Vert m^{\prime }(\alpha -\alpha
_{0})\right\Vert _{\mathcal{B}}>L\left\Vert \alpha -\alpha _{0}\right\Vert _{%
\mathcal{A}}$. In general, by $L\left\Vert \alpha -\alpha _{0}\right\Vert _{%
\mathcal{A}}^{r}\geq 0$, we have \textit{$\mathcal{N}^{\prime \prime \prime
}\subseteq $} $\mathcal{N}^{\prime }$ for $\mathcal{N}^{\prime }$ from
Assumption 1, so the rank condition is imposed by restricting attention to
the \textit{$\mathcal{N}$} of Theorem 2. Here the rank condition is still
important, since if it is not satisfied on some interesting set \textit{$%
\mathcal{N}^{\prime }$}, Theorem 2 cannot give local identification on an
interesting set \textit{$\mathcal{N}$}.

Theorem 2 forges a link between the curvature of $m(\alpha )$ as in
Assumption 2 and the identification set $\mathcal{N}.$ An example is a
scalar $\alpha $ and twice continuously differentiable $m(\alpha )$ with
bounded second derivative. Here Assumption 2 will be satisfied with $r=2,$ $%
L=\sup_{\alpha }\left\vert d^{2}m(\alpha )/d\alpha ^{2}\right\vert /2$, and $%
\mathit{\mathcal{N}}^{\prime \prime }$ equal to the real line, where $%
\left\vert \cdot \right\vert $ denotes the absolute value. Assumption 1 will
be satisfied with $\mathit{\mathcal{N}}^{\prime }$ equal to the real line as
long as $m^{\prime }=dm(\alpha _{0})/d\alpha $ is nonzero. Then $\mathcal{N}%
\mathit{^{\prime \prime \prime }}=\{\alpha :\left\vert \alpha -\alpha
_{0}\right\vert <L^{-1}\left\vert m^{\prime }\right\vert \}.$ Here $%
L^{-1}\left\vert m^{\prime }\right\vert $ is the minimum distance $\alpha $
must go from $\alpha _{0}$ before $m(\alpha )$ can "bend back" to zero. In
nonparametric models $\mathcal{N}^{\prime \prime \prime }$ will be an
analogous set.

When $r=1$ the set $\mathcal{N}^{\prime \prime \prime }$ will be a linear
cone with vertex at $\alpha _{0}$, which means that if $\alpha \in \mathcal{N%
}^{\prime \prime \prime }$ then so is $\lambda \alpha +(1-\lambda )\alpha
_{0}$ for $\lambda >0$. In general, $\mathcal{N}^{\prime \prime \prime }$ is
not convex, so it is not a convex cone. For $r>1$ the set $\mathit{\mathcal{N%
}^{\prime \prime \prime }}$ is not a cone although it is star shaped around $%
\alpha _{0}$, meaning that for any $\alpha \in \mathit{\mathcal{N}^{\prime
\prime \prime }}$ we have $\lambda \alpha +(1-\lambda )\alpha _{0}\in 
\mathit{\mathcal{N}^{\prime \prime \prime }}$ for $0<\lambda \leq 1.$

Also, if $r>1$ then for any $L>0$ and $1\leq r^{\prime }<r$ there is $\delta
>0$ such that 
\begin{equation*}
\mathit{\mathcal{N}}_{\delta }\cap \{\alpha :\left\Vert m^{\prime }(\alpha
-\alpha _{0})\right\Vert _{\mathcal{B}}>L\left\Vert \alpha -\alpha
_{0}\right\Vert _{\mathcal{A}}^{r^{\prime }}\}\subseteq \mathit{\mathcal{N}}%
_{\delta }\cap \mathit{\mathcal{N}^{\prime \prime \prime }.}
\end{equation*}%
In this sense $\alpha \in \mathit{\mathcal{N}^{\prime \prime \prime }}$ as
assumed in Theorem 2 is less restrictive the larger is $r$, i.e. the local
identification neighborhoods of Theorem 2 are "richer" the larger is $r$.

\subsection{Discussion of Assumptions 1 and 2}

Restricting the set of $\alpha $ to be smaller than an open ball can be
necessary for local identification in nonparametric models, as we now show
in an example. Suppose $\alpha =(\alpha _{1},\alpha _{2},...)$ is a sequence
of real numbers. Let $(p_{1},p_{2},...)$ be probabilities, $p_{j}>0,$ $%
\sum_{j=1}^{\infty }p_{j}=1.$ Let $f(x)$ be a twice continuously
differentiable function of a scalar $x$ that is bounded with bounded second
derivative. Suppose $f(x)=0$ if and only if $x\in \{0,1\}$ and $df(0)/dx=1$.
Let $m(\alpha )=(f(\alpha _{1}),f(\alpha _{2}),...)$ also be a sequence with 
$\left\Vert m(\alpha )\right\Vert _{\mathcal{B}}=\left( \sum_{j=1}^{\infty
}p_{j}f(\alpha _{j})^{2}\right) ^{1/2}$. Then for $\left\Vert \alpha
\right\Vert _{\mathcal{A}}=$ $\left( \sum_{j=1}^{\infty }p_{j}\alpha
_{j}^{4}\right) ^{1/4}$ the function $m(\alpha )$ will be Fr\'{e}chet
differentiable at $\alpha _{0}=0,$ with $m^{\prime }h=h$. A fourth moment
norm for $\alpha ,$ rather than a second moment norm, is needed to make $%
m(\alpha )$ Fr\'{e}chet differentiable under the second moment norm for $%
m(\alpha )$. Here the map $m^{\prime }$ is not onto, even though it is the
identity, because the norm on $\mathcal{A}$ is stronger than the norm on $%
\mathcal{B}$.

In this example the value $\alpha _{0}=0$ is not locally identified by the
equation $m(\alpha )=0$ on any open ball in the norm $\left\Vert \alpha
-\alpha _{0}\right\Vert _{\mathcal{A}}$. To show this result consider $%
\alpha ^{k}$ which has zeros in the first $k$ positions and a one everywhere
else, i.e., $\alpha ^{k}=(0,...,0,1,1,...).$ Then $m(\alpha ^{k})=0$ and for 
$\Delta ^{k}=\sum_{j=k+1}^{\infty }p_{j}\longrightarrow 0$ we have $%
\left\Vert \alpha ^{k}-\alpha _{0}\right\Vert _{\mathcal{A}}=\left(
\sum_{j=1}^{\infty }p_{j}[\alpha _{j}^{k}]^{4}\right) ^{1/4}=\left( \Delta
^{k}\right) ^{1/4}\longrightarrow 0.$ Thus, we have constructed a sequence
of $\alpha ^{k}$ not equal to $\alpha _{0}$ such that $m(\alpha ^{k})=0$ and 
$\left\Vert \alpha ^{k}-\alpha _{0}\right\Vert _{\mathcal{A}}\longrightarrow
0.$

We can easily describe the set $\mathcal{N}$ of Theorem 2 in this example,
on which $\alpha _{0}=0$ will be locally identified. By the second
derivative of $f$ being bounded, Assumption 2\textbf{\ }is satisfied with $%
\mathcal{N}^{\prime \prime }=\mathcal{A}$, $r=2,$ and $L=\sup_{a}|\partial
^{2}f(a)/\partial a^{2}|/2,$ where $L\geq 1$ by the fact that $f^{\prime
}(0)=1$ and $f(0)=f(1)=0$ (an expansion gives $0=f(1)=1+2^{-1}\partial ^{2}f(%
\bar{a})/\partial a^{2}$ for $0\leq \bar{a}\leq 1$)$.$ Then,%
\begin{equation*}
\mathcal{N=}\left\{ \alpha =(\alpha _{1},\alpha _{2},...):\left(
\sum_{j=1}^{\infty }p_{j}\alpha _{j}^{2}\right) ^{1/2}>L\left(
\sum_{j=1}^{\infty }p_{j}\alpha _{j}^{4}\right) ^{1/2}\right\} .
\end{equation*}%
The sequence $(\alpha ^{k})_{k=1}^{\infty }$ given above will not be
included in this set because $L\geq 1.$ A simple subset of $\mathcal{N}$ (on
which $\alpha _{0}$ is locally identified) is $\left\{ \alpha =(\alpha
_{1},\alpha _{2},...):\left\vert \alpha _{j}\right\vert
<L^{-1},(j=1,2,...)\right\} $.

It is important to note that Theorems 1 and 2 provide sufficient, and not
necessary, conditions for local identification. In fact, the conditions of
Theorems 1 and 2 are sufficient for%
\begin{equation}
\left\Vert m(\alpha )-m^{\prime }(\alpha -\alpha _{0})\right\Vert _{\mathcal{%
B}}\neq \left\Vert m^{\prime }(\alpha -\alpha _{0})\right\Vert _{\mathcal{B}}
\label{tang cone cond}
\end{equation}%
that implies $m(\alpha )\neq 0,$ to hold on $\mathcal{N}$. The set where (%
\ref{tang cone cond}) holds may be larger than the set $\mathcal{N}$ of
Theorem 1 or 2. We have focused on the set $\mathcal{N}$ of Theorem 1 or 2
because those conditions and the associated locally identified set $\mathcal{%
N}$ are relatively easy to interpret. See Appendix E for more identification
results related to (\ref{tang cone cond}).

Assumption 1 may not be needed for identification in nonlinear models,
although local identification is complicated in the absence of Assumption 1.
Conditions may involve nonzero higher order derivatives. Such results for
parametric models are discussed by, e.g., Sargan (1983). Here we focus on
models where Assumption 1 is satisfied.

\section{Local Identification in Hilbert Spaces}

\subsection{Full Rank Condition in Hilbert Spaces}

The restrictions imposed on $\alpha $ in Theorem 2 are not very transparent.
In Hilbert spaces it is possible to give more interpretable conditions based
on a lower bound for $\left\Vert m^{\prime }(\alpha -\alpha _{0})\right\Vert
_{\mathcal{B}}^{2}$. Let $\langle \cdot ,\cdot \rangle $ denote the inner
product for a Hilbert space.

\bigskip

\textsc{Assumption 3: }$(\mathcal{A},\left\Vert \cdot \right\Vert _{\mathcal{%
A}})$ \textit{and }$(\mathcal{B},\left\Vert \cdot \right\Vert _{\mathcal{B}%
}) $\textit{\ are separable Hilbert spaces and either }a) \textit{there is a
set }$\mathcal{N}^{\prime },$ \textit{an orthonormal basis }$\{\phi
_{1},\phi _{2},...\}\subseteq \mathcal{A},$\textit{\ and a bounded, positive
sequence }$(\mu _{1},\mu _{2},...)$ \textit{such that for all }$\alpha \in 
\mathcal{N}^{\prime },$%
\begin{equation*}
\left\Vert m^{\prime }(\alpha -\alpha _{0})\right\Vert _{\mathcal{B}%
}^{2}\geq \sum_{j=1}^{\infty }\mu _{j}^{2}\langle \alpha -\alpha _{0},\phi
_{j}\rangle ^{2};
\end{equation*}%
\textit{or b) }$m^{\prime }$\textit{\ is a compact linear operator with
positive singular values }$(\mu _{1},\mu _{2},...)$.

\bigskip

The hypothesis in b) that $m^{\prime }$\ is a compact operator is a mild one
when $m^{\prime }$ is a conditional expectation. Recall that an operator $m:%
\mathcal{A}\mapsto \mathcal{B}$ is compact if and only if it is continuous
and maps bounded sets in $\mathcal{A}$ into relatively compact sets in $%
\mathcal{B}$. Under very mild conditions, $m(\alpha )=E[\alpha (X)|W]$ is
compact: See Zimmer (1990, chapter 3), Kress (1999, section 2.4) and
Carrasco, Florens, and Renault (2007) for a variety of sufficient
conditions. When $m^{\prime }$ in b) is compact there is an orthonormal
basis $\{\phi _{j}:j=1,\ldots \}$ for $\mathcal{A}$ with 
\begin{equation*}
\left\Vert m^{\prime }(\alpha -\alpha _{0})\right\Vert _{\mathcal{B}%
}^{2}=\sum_{j=1}^{\infty }\mu _{j}^{2}\langle \alpha -\alpha _{0},\phi
_{j}\rangle ^{2},
\end{equation*}%
where $\mu _{j}^{2}$ are the eigenvalues and $\phi _{j}$\ the eigenfunctions
of the operator $m^{\prime \ast }m^{\prime }$, so that condition a) is
satisfied, where $m^{\prime \ast }$ denotes the adjoint of $m^{\prime }$.
The assumption that the singular values are all positive implies the rank
condition holds for $\mathcal{N}^{\prime }=\mathcal{A}$. Part a) differs
from part b) by imposing a lower bound on $\left\Vert m^{\prime }(\alpha
-\alpha _{0})\right\Vert _{\mathcal{B}}^{2}$ only over a subset $\mathcal{N}%
^{\prime }$ of $\mathcal{A}$ and by allowing the basis $\{\phi _{j}\}$\ to
be different from the eigenfunction basis of the operator $m^{\prime \ast
}m^{\prime }$. In principle this allows us to impose restrictions on $\alpha
-\alpha _{0}$, like boundedness and smoothness, which could help Assumption
3 a) to hold. For similar assumptions in estimation context, see, e.g., Chen
and Rei{\ss } (2011) and Chen and Pouzo (2012).

It turns out that there is a precise sense in which the rank condition is
satisfied for most data generating processes, if it is satisfied for one, in
the Hilbert space environment here. In this sense the rank condition turns
out to be generic. Let $\mathcal{A}$ and $\mathcal{B}$ be separable Hilbert
spaces, and $\mathcal{N}^{\prime }\subseteq \mathcal{A}$. Suppose that there
exists at least one compact linear operator: $\mathcal{K}:\mathcal{A}\mapsto 
\mathcal{B}$ which is injective , i.e. $\mathcal{K}\delta =0$ for $\delta
\in \mathcal{A}$ if and only if $\delta =0$. This is an infinite-dimensional
analog of the order condition, that for example rules out $\mathcal{B}$
having smaller finite dimension than $\mathcal{A}$ (e.g. having fewer
instruments than right-hand side endogenous variables in a linear regression
model). The operator $m^{\prime }:\mathcal{N}^{\prime }\mapsto \mathcal{B}$
is generated by the nature as follows:

\begin{itemize}
\item[1.] The nature selects a countable orthonormal basis $\{\phi _{j}\}$
of cardinality $N\leq \infty $ in $\mathcal{A}$ and an orthonormal set $%
\{\varphi _{j}\}$ of equal cardinality in $\mathcal{B}$.

\item[2.] The nature samples a bounded sequence of real numbers $\{\lambda
_{j}\}$ according to a probability measure $\eta $ whose each marginal is
dominated by the Lebesgue measure on $\mathbb{R}$, namely $\text{Leb}(A)=0$
implies $\eta (\{\lambda _{j}\in A\})=0$ for any measurable $A\subset 
\mathbb{R}$ for each $j$.
\end{itemize}

Then the nature sets, for some scalar number $\kappa >0$, and every $\delta
\in \mathcal{N}^{\prime }$ 
\begin{equation}
m^{\prime }\delta =\kappa (\sum_{j=0}^{N}\lambda _{j}\langle \phi
_{j},\delta \rangle \varphi _{j}).
\end{equation}%
This operator is properly defined on $\mathcal{N}^{\prime }:=\{\delta \in 
\mathcal{A}:m^{\prime }\delta \in \mathcal{B}\}$. \bigskip

\textsc{Lemma 3} \textit{(1) In the absence of further restrictions on $%
m^{\prime }$, the algorithms obeying conditions 1 and 2 exist. (2) If $%
m^{\prime }$ is generated by any algorithm that obeys conditions 1 and 2,
then probability that $m^{\prime }$ is not injective over }$\mathcal{N}%
^{\prime }$\textit{\ is zero, namely} $\text{Pr}_{\eta }\{\exists \delta \in 
\mathcal{N}^{\prime }:\delta \neq 0\text{ and }m^{\prime }\delta =0\}=0.$ 
\textit{Moreover, Assumption 3 holds with $\mu _{j}=|\kappa \lambda _{j}|$
with probability one under $\eta $. }

\bigskip

In Appendix B we provide examples for the case $\mathcal{A}=\mathcal{B}%
=L^{2}[0,1]$ that highlight the range of algorithms permitted by conditions
1 and 2 above, including cases where various restrictions on $m^{\prime }$
are imposed: boundedness, compactness, weak positivity, and density
restrictions. Genericity arguments use the idea of randomization, and are
often employed in economic theory, functional analysis, and probability
theory, see, e.g., Anderson and Zame (2000), Marcus and Pisier (1981),
Ledoux and Talagrand (2011). Andrews (2011) previously used a related notion
of genericity, called prevalence within bounded sets, to argue that rich
classes of operators induced by densities in nonparametric IV are $L^{2}$%
-complete. Though inspired in part by Andrews (2011), the simple result
above uses a somewhat different notion of genericity than prevalence\textbf{.%
}\footnote{%
Informally speaking, prevalence requires that it should be possible to
construct a randomization device such that \textit{all} finite-dimensional
distributions for $\lambda _{j}$'s are absolutely continuous, i.e., the
distribution of $(\lambda _{j_{1}},...,\lambda _{j_{k}})$ needs to be
continuous with respect to the Lebesgue measure on $\mathbb{R}^{k}$, for any 
$(j_{1},...,j_{k})\subset \{1,2,...,N\}$, and any $k\in \{1,2,...\}$. The
notion that we use requires \textit{only} that the \textit{one-dimensional}
marginal distributions for $\lambda _{j}$ are absolutely continuous for any $%
j$. The distinction is actually important to cover cases, where perfect
dependence between some $\lambda _{j}$'s may be required to maintain
conditions imposed on the operator, such as, for example, the kernel of the
operator being a conditional density. See Appendix B for details.} We also
note that while this construction implies identification with probability
one, it does not regulate in any way the strength of identification, and
hence has no bearing on the choice of an inferential method.%

\subsection{Local Identification in Hilbert Spaces}

In what follows let $b_{j}=\langle \alpha -\alpha _{0},\phi _{j}\rangle
,j=1,2,...$ denote the Fourier coefficients for $\alpha -\alpha _{0}$, so
that $\alpha =\alpha _{0}+\sum_{j=1}^{\infty }b_{j}\phi _{j}$. Under
Assumptions 2 and 3 we can characterize an identified set in terms of the
Fourier coefficients.

\bigskip

\textsc{Theorem 4:}\textit{\ If Assumptions 2 and 3 are satisfied} \textit{%
then }$\alpha _{0}$ \textit{is locally identified on $\mathcal{N=N}^{\prime
\prime }\cap \mathcal{N}^{\prime \prime \prime }$}\textbf{, }\textit{where }%
\begin{equation*}
\mathit{\mathcal{N}^{\prime \prime \prime }=}\{\alpha =\alpha
_{0}+\sum_{j=1}^{\infty }b_{j}\phi _{j}:\sum_{j=1}^{\infty }\mu
_{j}^{2}b_{j}^{2}>L^{2}(\sum_{j=1}^{\infty }b_{j}^{2})^{r}\}.
\end{equation*}

\bigskip

When $r=1$ it is necessary for $\alpha \in \mathit{\mathcal{N}^{\prime
\prime \prime }}$ that the Fourier coefficients $b_{j}$ where $\mu _{j}^{2}$
is small not be too large relative to the Fourier coefficients where $\mu
_{j}^{2}$ is large. In particular, when $r=1$ any $\alpha \neq \alpha _{0}$
with $b_{j}=0$ for all $j$ with $\mu _{j}>L$ will not be an element of $%
\mathit{\mathcal{N}^{\prime \prime \prime }}.$ When $r>1$ we can use the H%
\"{o}lder inequality to obtain a sufficient condition for $\alpha \in 
\mathit{\mathcal{N}^{\prime \prime \prime }}$ that is easier to interpret.

\bigskip

\textsc{Corollary 5:}\textit{\ If Assumptions 2 and 3 are satisfied, with }$%
L>0,r>1,$ \textit{then }$\alpha _{0}$\textit{\ is locally identified on }$%
\mathcal{N}=\mathcal{N}^{\prime \prime }\mathcal{\cap N}^{\prime \prime
\prime }$ $\mathit{where}$ $\mathcal{N}^{\prime \prime \prime }=\{\alpha
=\alpha _{0}+\sum_{j=1}^{\infty }b_{j}\phi _{j}:\sum_{j=1}^{\infty }\mu
_{j}^{-2/(r-1)}b_{j}^{2}<L^{-2/(r-1)}\}.$

\bigskip

For $\alpha $ to be in the $\mathcal{N}^{\prime \prime \prime }$ of
Corollary 5 the Fourier coefficients $b_{j}$ must vanish faster than $\mu
_{j}^{1/(r-1)}$ as $j$ grows. In particular, a sufficient condition for $%
\alpha \in \mathcal{N}^{\prime \prime \prime }$ is that $\left\vert
b_{j}\right\vert <\left( \mu _{j}/L\right) ^{1/(r-1)}c_{j}$ for any positive
sequence $c_{j}$ with $\sum_{j=1}^{\infty }c_{j}^{2}=1$. These bounds on $%
b_{j}$ correspond to a hyperrectangle while the $\mathcal{N}^{\prime \prime
\prime }$ in Corollary 5 corresponds to an ellipsoid. The bounds on $b_{j}$
shrink as $L$ increases, corresponding to a smaller local identification set
when more nonlinearity is allowed. Also, it is well known that, at least in
certain environments, imposing bounds on Fourier coefficients corresponds to
imposing smoothness conditions, like existence of derivatives; see for
example Kress (Chapter 8, 1999). In that sense the identification set in
Corollary 5 imposes smoothness conditions on the deviations of $\alpha $
from the truth $\alpha _{0}$.

The bound imposed in $\mathcal{N}^{\prime \prime \prime }$ of Corollary 5 is
a ``source condition" under Assumption 3 b) and is similar to conditions
used by Florens, Johannes and Van Bellegem (2011) and others. Under
Assumption 3 a) it is similar to norms in generalized Hilbert scales, for
example, see Engl, Hanke, and Neubauer (1996) and Chen and Rei{\ss } (2011).
Our Assumption 3 a) or 3 b) are imposed on deviations $\alpha -\alpha _{0}$,
while the above references all impose on true function $\alpha _{0}$\ itself
as well as on the parameter space hence on the deviations.

\subsection{A Quantile IV Example}

To illustrate the results of this Section we consider an endogenous quantile
example where $0<\tau <1$ is a scalar,%
\begin{equation*}
\rho (Y,X,\alpha )=1(Y\leq \alpha (X))-\tau ,
\end{equation*}%
$\mathcal{A}=\{\alpha (\cdot ):\text{E}[\alpha (X)^{2}]<\infty \}$, and $%
\mathcal{B}=\{a(\cdot ):\text{E}[a(W)^{2}]<\infty \},$ with the usual
Hilbert spaces of mean-square integrable random variables. Here we have 
\begin{equation*}
m(\alpha )=\text{E}[1(Y\leq \alpha (X))|W]-\tau .
\end{equation*}%
Let $f_{Y}(y|X,W)$ denote the conditional probability density function (pdf)
of $Y$ given $X$ and $W$, \textit{\ }$f_{X}(x|W)$\textbf{\ }the conditional
pdf of $X$\ given $W$, and $f(x)$\ the marginal pdf of $X$.

\bigskip

\textsc{Theorem 6: }\textit{If }$f_{Y}(y|X,W)$\textit{\ is continuously
differentiable in }$y$ \textit{with }$\left\vert df_{Y}(y|X,W)/dy\right\vert
\leq L_{1},$\textit{\ }$f_{X}(x|W)\leq L_{2}f(x),$ \textit{and }$m^{\prime
}h=\text{E}[f_{Y}(\alpha _{0}(X)|X,W)h(X)|W]$\textit{\ satisfies Assumption
3, then }$\alpha _{0}$\textit{\ is locally identified on }%
\begin{equation*}
\mathcal{N}=\{\alpha =\alpha _{0}+\sum_{j=1}^{\infty }b_{j}\phi _{j}\in 
\mathcal{A}:\sum_{j=1}^{\infty }b_{j}^{2}/\mu _{j}^{2}<(L_{1}L_{2})^{-2}\}.
\end{equation*}

\bigskip

This result gives a precise link between a neighborhood on which $\alpha
_{0} $ is locally identified and the bounds $L_{1}$ and $L_{2}$. Assumption
3 b) will hold under primitive conditions for $m^{\prime }$\ to be complete,
that are given by Chernozhukov, Imbens, and Newey (2007). Theorem 6 corrects
Theorem 3.2 of Chernozhukov, Imbens, and Newey (2007) by adding the bound on 
$\sum_{j=1}^{\infty }b_{j}^{2}/\mu _{j}^{2}.$ It also gives primitive
conditions for local identification for general $X$ while Chernozhukov and
Hansen (2005) only gave primitive conditions for identification when $X$ is
discrete. Horowitz and Lee (2007) impose analogous conditions in their paper
on convergence rates of nonparametric endogenous quantile estimators but
assumed identification.

\section{Semiparametric Models}

\subsection{Identification Results}

In this section, we consider local identification in possibly nonlinear
semiparametric models, where $\alpha $ can be decomposed into a $p\times 1$
dimensional parameter vector $\beta $ and nonparametric component $g$, so
that $\alpha =(\beta ,g)$. Let $|\cdot |$ denote the Euclidean norm for $%
\beta $ and assume $g\in \mathcal{G}$ where $\mathcal{G}$ is a Banach space
with norm $\left\Vert \cdot \right\Vert _{\mathcal{G}},$ such as a Hilbert
space. We focus here on a conditional moment restriction model 
\begin{equation*}
\text{E}[\rho (Y,X,\beta _{0},g_{0})|W]=0,
\end{equation*}%
where $\rho (y,x,\beta ,g)$ is a $J\times 1$ vector of residuals. Here $%
m(\alpha )=\text{E}[\rho (Y,X,\beta ,g)|W]$ will be considered as an element
of the Hilbert space $\mathcal{B}$ of $J\times 1$ random vectors with inner
product%
\begin{equation*}
\left\langle a,b\right\rangle =\text{E}[a(W)^{T}b(W)].
\end{equation*}

The differential $m^{\prime }(\alpha -\alpha _{0})$ can be expressed as 
\begin{equation*}
m^{\prime }(\alpha -\alpha _{0})=m_{\beta }^{\prime }(\beta -\beta
_{0})+m_{g}^{\prime }(g-g_{0}),
\end{equation*}%
where $m_{\beta }^{\prime }$ is the derivative of $m(\beta ,g_{0})=\text{E}%
[\rho (Y,X,\beta ,g_{0})|W]$ with respect to $\beta $ at $\beta _{0}$ and ${m%
}_{g}^{\prime }$ is the G\^{a}teaux derivative of $m(\beta _{0},g)$ with
respect to $g$ at $g_{0}.$ To give conditions for local identification of $%
\beta _{0}$ in the presence of the nonparametric component $g$ it is helpful
to partial out $g$. Let $\overline{\mathcal{M}}$ be the closure of the
linear span $\mathcal{M}$ of $m_{g}^{\prime }(g-g_{0})$ for $g\in \mathcal{N}%
_{g}^{\prime }$ where $\mathcal{N}_{g}^{\prime }$ will be specified below.
In general $\overline{\mathcal{M}}\mathcal{\neq M}$ because the linear
operator $m_{g}^{\prime }$ need not have closed range (like $m^{\prime }$
onto, a closed range would also imply a continuous inverse, by the Banach
inverse theorem). For the $k^{th}$ unit vector $e_{k},(k=1,...,p),$ let 
\begin{equation*}
\zeta _{k}^{\ast }=\arg \min_{\zeta \in \overline{\mathcal{M}}}\text{E}%
[\{m_{\beta }^{\prime }(W)e_{k}-\zeta (W)\}^{T}\{m_{\beta }^{\prime
}(W)e_{k}-\zeta (W)\}],
\end{equation*}%
which exists and is unique by standard Hilbert space results; e.g. see
Luenberger (1969). Define $\Pi $ to be the $p\times p$ matrix with 
\begin{equation*}
\Pi _{jk}:=\text{E}\left[ \left\{ m_{\beta }^{\prime }(W)e_{j}-\zeta
_{j}^{\ast }(W)\right\} ^{T}\left\{ m_{\beta }^{\prime }(W)e_{k}-\zeta
_{k}^{\ast }(W)\right\} \right] ,\text{\quad }(j,k=1,...,p).
\end{equation*}%
The following condition is important for local identification of $\beta _{0}$%
.

\bigskip

\textsc{Assumption 4: }$m^{\prime }:\mathbb{R}^{p}\times \mathcal{N}%
_{g}^{\prime }\longrightarrow \mathcal{B}$ \textit{is linear and bounded,
and }$\Pi $ \textit{is nonsingular.}


\bigskip

This assumption is similar to those first used by Chamberlain (1992) to
establish the possibility of estimating parametric components at root-$n$
rate in semi-parametric moment condition problems; see also Ai and Chen
(2003) and Chen and Pouzo (2009). In the local identification analysis
considered here it leads to local identification of $\beta _{0}$ without
identification of $g$ when $m(\beta _{0},g)$ is linear in $g$. It allows us
to separate conditions for identification of $\beta _{0}$ from conditions
for identification of $g.$ Note that the parameter $\beta $ may be
identified even when $\Pi $ is singular, but that case is more complicated,
as discussed at the end of Section 2, and we do not analyze this case.


The following condition controls the behavior of the derivative with respect
to $\beta $:

\bigskip

\textsc{Assumption 5: }\textit{For every }$\varepsilon >0$ \textit{there is
a neighborhood }$B$\textit{\ of }$\beta _{0}$\textit{\ and a set }$\mathcal{N%
}_{g}^{\beta }$\textit{\ such that for all }$g\in \mathcal{N}_{g}^{\beta }$ 
\textit{with probability one }$\text{E}[\rho (Y,X,\beta ,g)|W]$ \textit{is
continuously differentiable in }$\beta $ \textit{on }$B$\textit{\ and}%
\begin{equation*}
\sup_{g\in \mathcal{N}_{g}^{\beta }}\sqrt{\text{E}[\sup_{\beta \in
B}\left\vert \partial \text{E}[\rho (Y,X,\beta ,g)|W]/\partial \beta
-\partial \text{E}[\rho (Y,X,\beta _{0},g_{0})|W]/\partial \beta \right\vert
^{2}]}<\varepsilon .
\end{equation*}

\bigskip

It turns out that Assumptions 4 and 5 will be sufficient for local
identification of $\beta _{0}$ when $m(\beta _{0},g)$ is linear in $g,$ i.e.
for $m(\beta ,g)=0$ to imply $\beta =\beta _{0}$ when $(\beta ,g)$ is in
some neighborhood of $(\beta _{0},g_{0}).$ This works because Assumption 4
partials out the effect of unknown $g$ on local identification of $\beta
_{0} $.

\bigskip

\textsc{Theorem 7: }\textit{If Assumptions 4 and 5 are satisfied and }$%
m(\beta _{0},g)$ \textit{is linear in }$g$ \textit{then there is an }$%
\varepsilon >0$ \textit{such that for }$B$ \textit{and }$\mathcal{N}%
_{g}^{\beta }$ \textit{from Assumption 5 and }$\mathcal{N}_{g}^{\prime }$%
\textit{\ from Assumption 4, }$\beta _{0}$\textit{\ is locally identified
for }$\mathcal{N}=B\times (\mathcal{N}_{g}^{\prime }\cap \mathcal{N}%
_{g}^{\beta }).$ \textit{If, in addition, Assumption 1 is satisfied} \textit{%
for }$m_{g}^{\prime }$ \textit{and }$\mathcal{N}_{g}^{\prime }\cap \mathcal{N%
}_{g}^{\beta }$\textit{\ replacing }$m^{\prime }$\textit{\ and }$\mathcal{N}%
^{\prime }$\textit{\ then }$\alpha _{0}=(\beta _{0},g_{0})$\textit{\ is
locally identified for }$\mathcal{N}$.

\bigskip

This result is more general than Florens, Johannes, and Van Bellegem (2012)
and Santos (2011) since it allows for nonlinearities in $\beta $, and
dependence on $g$ of the partial derivatives $\partial \text{E}[\rho
(Y,X,\beta ,g)|W]/\partial \beta $. When the partial derivatives $\partial 
\text{E}[\rho (Y,X,\beta ,g)|W]/\partial \beta $ do not depend on $g$, then
Assumption 5 could be satisfied with $\mathcal{N}_{g}^{\prime }=\mathcal{G}$%
, and Theorem 7 could then imply local identification of $\beta _{0}$ in
some neighborhood of $\beta _{0}$ only.

For semiparametric models that are nonlinear in $g$ we can give local
identification results based on Theorem 2 or the more specific conditions of
Theorem 4 and Corollary 5. For brevity we give just a result based on
Theorem 2.

\bigskip

\textsc{Theorem 8: }\textit{If Assumptions 4 and 5 are satisfied and }$%
m(\beta _{0},g)$ \textit{satisfies Assumption 2 with $\mathcal{N}^{\prime
\prime }=\mathcal{N}_{g}^{\prime \prime }$, then there is an }$\varepsilon
>0 $ \textit{such that for }$B$ \textit{and }$\mathcal{N}_{g}^{\beta }$ 
\textit{from Assumption 5, }$\mathcal{N}_{g}^{\prime }$\textit{\ from
Assumption 4, and }%
\begin{equation*}
\mathit{\mathcal{N}_{g}^{\prime \prime \prime }}=\{g:\left\Vert
m_{g}^{\prime }(g-g_{0})\right\Vert _{\mathcal{B}}>\varepsilon
^{-1}L\left\Vert g-g_{0}\right\Vert _{\mathcal{A}}^{r}\}
\end{equation*}%
\textit{it is the case that }$\alpha _{0}=(\beta _{0},g_{0})$\textit{\ is
locally identified for }$\mathcal{N}=B\times (\mathcal{N}_{g}^{\beta }\cap 
\mathcal{N}_{g}^{\prime }\cap \mathcal{N}_{g}^{\prime \prime }\cap \mathcal{N%
}_{g}^{\prime \prime \prime }).$

\subsection{A Single Index IV Example}

Econometric applications often have too many covariates for fully
nonparametric estimation to be practical, i.e. they suffer from the curse of
dimensionality. Econometric practice thus motivates interest in models with
reduced dimension. An important such model is the single index model. Here
we consider a single index model with endogeneity, given by 
\begin{equation}
Y=g_{0}(X_{1}+X_{2}^{T}\beta _{0})+U,\text{\quad }\text{E}[U|W]=0,
\label{sing index}
\end{equation}%
where $\beta _{0}$ is a vector of unknown parameters, $g_{0}(\cdot )$ is an
unknown function, and $W$ are instrumental variables. Here the nonparametric
part is just one dimensional rather than having the same dimension as $X$.
This model is nonlinear in Euclidean parameters, and so is an example where
our results apply. Our results add to the literature on dimension reduction
with endogeneity, by showing how identification of an index model requires
fewer instrumental variables than a fully nonparametric IV model. We could
generalize the results to multiple indices but focus on a single index for
simplicity.

The location and scale of the parametric part are not identified separately
from $g_{0}$, and hence, we normalize the constant to zero and the
coefficient of $X_{1}$ to $1$. Here 
\begin{equation*}
m(\alpha )(W)=\text{E}[Y-g(X_{1}+X_{2}^{T}\beta )|W].
\end{equation*}%
Let $V=X_{1}+X_{2}^{T}\beta _{0}$ and for differentiable $g_{0}(V)$ let 
\begin{equation*}
m_{\beta }^{\prime }=-\text{E}[g_{0}^{\prime }(V)X_{2}^{T}|W].
\end{equation*}%
Let $\zeta _{j}^{\ast }$ denote the projection of $m_{\beta }^{\prime
}e_{j}=-\text{E}[g_{0}^{\prime }(V)X_{2j}|W]$ on the mean-square closure of
the set $\{\text{E}[h(V)|W]:\text{E}[h(V)^{2}]<\infty \}$ and $\Pi $ the
matrix with $\Pi _{jk}=\text{E}[(m_{\beta }^{\prime }e_{j}-\zeta _{j}^{\ast
})(m_{\beta }^{\prime }e_{k}-\zeta _{k}^{\ast })].$

\bigskip

\textsc{Theorem 9:} \textit{Consider the model of equation (\ref{sing index}%
). If a) }$g_{0}(V)$\textit{\ is continuously differentiable with bounded
derivative }$g_{0}^{\prime }(V)$\textit{\ satisfying }$|g_{0}^{\prime }(%
\tilde{V})-g_{0}^{\prime }(V)|\leq C_{g}|\tilde{V}-V|$\textit{\ for some }$%
C_{g}>0,$ \textit{b) }$\text{E}[\left\vert X_{2}\right\vert ^{4}]<\infty $%
\textit{, and c) }$\Pi $\textit{\ is nonsingular, then there is a
neighborhood }$B$\textit{\ of }$\beta _{0}$\textit{\ and }$\delta >0$\textit{%
\ such that for} 
\begin{equation*}
\mathcal{N}_{g}^{\delta }=\{g:g(v)\text{ \textit{is continuously
differentiable and} }\sup_{v}|g^{\prime }(v)-g_{0}^{\prime }(v)|\leq \delta
\}
\end{equation*}%
$\beta _{0}$ \textit{is locally identified for }$\mathcal{N}=B\times 
\mathcal{N}_{g}^{\delta }.$ \textit{Furthermore, if there is }$\mathcal{N}%
_{g}^{\prime }$\textit{\ such that }$\text{E}[g(V)-g_{0}(V)|W]$\textit{\ is
bounded complete on the set }$\{g(V)-g_{0}(V):g\in \mathcal{N}_{g}^{\prime
}\}$\textit{\ then }$(\beta _{0},g_{0})$\textit{\ is locally identified for }%
$\mathcal{N}=B\times (\mathcal{N}_{g}^{\delta }\cap \mathcal{N}_{g}^{\prime
}).$

\bigskip

Since this model includes as a special case the linear simultaneous
equations model the usual rank and order conditions are still necessary for $%
\Pi $ to be nonsingular for all possible models, and hence are necessary for
identification. Relative to the linear nonparametric IV model in Newey and
Powell (2003) the index structure lowers the requirements for identification
by requiring that $m_{g}^{\prime }h=-\text{E}[h(V)|W]$ be complete on $%
\mathcal{N}_{g}^{\prime }$ rather than completeness of the conditional
expectation of functions of $X$ given $W$. For example, it may be possible
to identify $\beta _{0}$ and $g_{0}$ with only two instrumental variables,
one of which is used to identify $g_{0}$ and functions of the other being
used to identify $\beta _{0}$.

To further explain we can give more primitive conditions for nonsingularity
of $\Pi $. The following result gives a necessary condition for $\Pi $ to be
nonzero (and hence nonsingular) as well as a sufficient condition for
nonsingularity of $\Pi $.

\bigskip

\textsc{Theorem 10:} \textit{Consider the model of (\ref{sing index}). If }$%
\Pi $\textit{\ is nonsingular then the conditional distribution of }$W$%
\textit{\ given }$V$ \textit{is not complete. Also, if there is a measurable
function }$T(W)$ \textit{such that the conditional distribution of }$V$ 
\textit{given }$W$ \textit{depends only on }$T(W)$ \textit{and for every }$%
p\times 1\mathit{\ }$\textit{vector }$\lambda \neq 0,$\textit{\ }$\text{E}%
[g_{0}^{\prime }(V)\lambda ^{T}X_{2}|W]$\textit{\ is not measurable with
respect to }$T(W)$, \textit{then }$\Pi $ \textit{is nonsingular.}

\bigskip

To explain the conditions of this result note that if there is only one
variable in $W$ then the completeness condition (of $W$\ given $V$) can hold
and hence $\Pi $ can be singular. If there is more than one variable in $W$
then generally completeness (of $W$\ given $V$) will not hold, because
completeness would be like identifying a function of more than one variable
(i.e. $W$) with one instrument (i.e. $V$). If $W$ and $V$ are joint Gaussian
and $V$ and $W$ are correlated then completeness holds (and hence $\Pi $ is
singular) when $W$ is one dimensional but not otherwise. In this sense
having more than one instrument in $W$ is a necessary condition for
nonsingularity of $\Pi $. Intuitively, one instrument is needed for
identification of the one dimensional function $g_{0}(V)$ so that more than
one instrument is needed for identification of $\beta $.

The sufficient condition for nonsingularity of $\Pi $ is stronger than
noncompleteness. It is essentially an exclusion restriction, where $\text{E}%
[g_{0}^{\prime }(V)X_{2}|W]$ depends on $W$ in a different way than the
conditional distribution of $V$ depends on $W$. This condition can be shown
to hold if $W$ and $V$ are Gaussian, $W$ is two dimensional, and $\text{E}%
[g_{0}^{\prime }(V)X_{2}|W]$ depends on all of $W$.

\section{Semiparametric CCAPM}

Consumption capital asset pricing models (CCAPM) provide interesting
examples of nonparametric and semiparametric moment restrictions; see
Gallant and Tauchen (1989), Newey and Powell (1988), Hansen, Heaton, Lee,
and Roussanov (2007), Chen and Ludvigson (2009), and others. In this
section, we apply our general theorems to develop new results on
identification of a particular semiparametric specification of marginal
utility of consumption. Our results could easily be extended to other
specifications, and so are of independent interest.

To describe the model let $C_{t}$ denote consumption level at time $t$ and $%
c_{t}\equiv C_{t}/C_{t-1}$ be consumption growth. Suppose that the marginal
utility of consumption at time $t$ is given by 
\begin{equation*}
MU_{t}=C_{t}^{-\gamma _{0}}g_{0}(C_{t}/C_{t-1})=C_{t}^{-\gamma
_{0}}g_{0}(c_{t}),
\end{equation*}%
where $g_{0}(c)$ is an unknown positive function. For this model the
intertemporal marginal rate of substitution is%
\begin{equation*}
\delta _{0}MU_{t+1}/MU_{t}=\delta _{0}c_{t+1}^{-\gamma
_{0}}g_{0}(c_{t+1})/g_{0}(c_{t}),
\end{equation*}%
where $0<\delta _{0}\leq 1$ is the rate of time preference. Let $%
R_{t+1}=(R_{t+1,1},...,R_{t+1,J})^{T}$ be a $J\times 1$ vector of gross
asset returns. A semiparametric CCAPM equation is then given by%
\begin{equation}
\text{E}[R_{t+1}\delta _{0}c_{t+1}^{-\gamma
_{0}}\{g_{0}(c_{t+1})/g_{0}(c_{t})\}|W_{t}]=e,  \label{pri eqn}
\end{equation}%
where $e$ is a $J\times 1$ vector of ones, and $W_{t}\equiv
(Z_{t}^{T},c_{t})^{T}$ is a vector of random variables observed by the agent
at time $t$, with $Z_{t}$ not a measurable function of $c_{t}$. This
corresponds to an external habit formation model with only one lag, a
special case of Chen and Ludvigson (2009). As emphasized in Cochrane (2005),
habit formation models can help explain the high risk premia embedded in
asset prices. We focus here on consumption growth $c_{t}=C_{t}/C_{t-1}$ to
circumvent the potential nonstationarity of the level of consumption, see
Hall (1978), as has long been done in this literature, e.g. Hansen and
Singleton (1982).

From economic theory it is known that under complete markets there is a
unique intertemporal marginal rate of substitution that solves equation (\ref%
{pri eqn}), when $R_{t}$ is allowed to vary over all possible vectors of
asset returns. Of course that does not guarantee a unique solution for a
fixed vector of returns $R_{t}$. Note though that the semiparametric model
does impose restrictions on the marginal rate of substitution that should be
helpful for identification. We show how these restrictions lead to local
identification of this model via the results of Section 4.

This model can be formulated as a semiparametric conditional moment
restriction by letting $Y=(R_{t+1}^{T},c_{t+1},c_{t})^{T},$ $\beta =(\delta
,\gamma )^{T}$, $W=W_{t}=(Z_{t}^{T},c_{t})^{T},$ and 
\begin{equation}
\rho (Y,\beta ,g)=R_{t+1}\delta c_{t+1}^{-\gamma }g(c_{t+1})-g(c_{t})e.
\label{capm resid}
\end{equation}%
Then multiplying equation (\ref{pri eqn}) through by $g_{0}(c_{t})$ gives
the conditional moment restriction $\text{E}[\rho (Y,\beta _{0},g_{0})|W]=0$%
. Let $A_{t}=R_{t+1}\delta _{0}c_{t+1}^{-\gamma _{0}}$. The nonparametric
rank condition (Assumption 1 for $g$\textbf{)} will be uniqueness, up to
scale, of the solution $g_{0}$ of 
\begin{equation}
\text{E}[A_{t}g(c_{t+1})|W_{t}]=g(c_{t})e.  \label{2nd kind eq}
\end{equation}%
This equation differs from the linear nonparametric IV restriction where the
function $g_{0}(X)$ would solve $\text{E}[Y|W]=\text{E}[g(X)|W]$. That
equation is an integral equation of the first kind while equation (\ref{2nd
kind eq}) is a homogeneous integral equation of the second kind. The rank
condition for this second kind equation is that the null space of the
operator $\text{E}[A_{t}g(c_{t+1})|W_{t}]-g(c_{t})e$ is one-dimensional,
which is different than the completeness condition for first kind equations.
This example illustrates that the rank condition of Assumption 1 need not be
equivalent to completeness of a conditional expectation. Escanciano and
Hoderlein (2010) and Lewbel, Linton, and Srisuma (2012) have previously
shown how homogenous integral equations of the second kind arise in CCAPM
models, though their models and identification results are different than
those given here, as further discussed below.

Let $X_{t}=(1/\delta _{0},-\ln (c_{t+1}))^{T}$. Then differentiating inside
the integral, as allowed under regularity conditions given below, and
applying the Gateaux derivative calculation gives%
\begin{equation*}
m_{\beta }^{\prime }(W)=\text{E}[A_{t}g_{0}(c_{t+1})X_{t}^{T}|W_{t}],\ \
m_{g}^{\prime }g=\text{E}[A_{t}g(c_{t+1})|W_{t}]-g(c_{t})e.
\end{equation*}%
When $\text{E}[A_{t}g(c_{t+1})|W_{t}]$ is a compact operator, as holds under
conditions described below, it follows from the theory of integral equations
of the second kind (e.g. Kress, 1999, Theorem 3.2) that the set of
nonparametric directions $\mathcal{M}$ will be closed, i.e. 
\begin{equation*}
\overline{\mathcal{M}}=\mathcal{M=\{}\text{E}%
[A_{t}g(c_{t+1})|W_{t}]-g(c_{t})e:\left\Vert g\right\Vert _{\mathcal{G}%
}<\infty \},
\end{equation*}%
where we will specify $\left\Vert g\right\Vert _{\mathcal{G}}$ below. Let $%
\Pi $ be the two-dimensional second moment matrix $\Pi $ of the residuals
from the projection of each column of $m_{\beta }^{\prime }$ on $\overline{%
\mathcal{M}},$ as described in Section 4. Then nonsingularity of $\Pi $
leads to local identification of $\beta _{0}$ via Theorem 7.

To give a precise result let $\Delta $ be any finite positive number, 
\begin{eqnarray*}
D_{t} &=&(1+\left\vert R_{t+1}\right\vert )[2+|\ln
(c_{t+1})|^{2}]\sup_{\gamma \in \lbrack \gamma _{0}-\Delta _{,}\gamma
_{0}+\Delta ]}c_{t+1}^{-\gamma }, \\
\mathcal{G} &=&\left\{ g:\left\Vert g\right\Vert _{\mathcal{G}}\equiv \sqrt{%
\text{E}[E[D_{t}^{2}|W_{t}]g(c_{t+1})^{2}]}<\infty \right\} .
\end{eqnarray*}%
The following assumption imposes some regularity conditions.

\bigskip

\textsc{Assumption 6:} $(R_{t}^{T},c_{t},Z_{t}^{T})$ \textit{is strictly
stationary, }$\text{E}[D_{t}^{2}]<\infty $; $0<\delta _{0}\leq 1$, $%
\left\Vert g_{0}\right\Vert _{\mathcal{G}}<\infty $.\textit{\ }

\bigskip

The following result applies Theorem 7 to this CCAPM.

\bigskip

\textsc{Theorem 11:} \textit{Consider equation (\ref{2nd kind eq}).\ Suppose
that Assumption 6 is satisfied. Then the linear mapping }$m^{\prime }:%
\mathbb{R}^{2}\times \mathcal{G}\longrightarrow \mathcal{B}$ \textit{is
bounded and if }$\Pi $ \textit{is nonsingular there is a neighborhood }$B$%
\textit{\ of }$\beta _{0}$\textit{\ and }$\varepsilon >0$\textit{\ such that
for} $\mathcal{N}_{g}^{\beta }=\{g:\left\Vert g-g_{0}\right\Vert _{\mathcal{G%
}}<\varepsilon \},$ $\beta _{0}$ \textit{is locally identified for }$%
\mathcal{N}=B\times \mathcal{N}_{g}^{\beta }.$ \textit{If in addition }$%
m_{g}^{\prime }(g-g_{0})\neq 0$\textit{\ for all }$g\neq g_{0}$\textit{\ and 
}$g\in \mathcal{N}_{g}^{\beta }$\textit{\ then }$(\beta _{0},g_{0})$ \textit{%
is locally identified for }$\mathcal{N}=B\times \mathcal{N}_{g}^{\beta }$.

\bigskip

Primitive conditions for nonsingularity of $\Pi $ and for $m_{g}^{\prime
}(g-g_{0})\neq 0$ when\textit{\ }$g\neq g_{0}$\textit{\ }are needed to make
this result interesting. It turns out that some completeness conditions
suffice, as shown by the following result. Let $\tilde{W}%
_{t}=(w(Z_{t}),c_{t})$ for some measurable function $w(Z_{t})$ of $Z_{t}$,
and $f_{c,\tilde{W}}(c,\tilde{w})$ denote the joint pdf of $(c_{t+1},\tilde{W%
}_{t})$, $f_{c}(c)$ and $f_{\tilde{W}}(\tilde{w})$\ the marginal pdfs of $%
c_{t+1}$ and $\tilde{W}_{t}$\ respectively.

\bigskip

\textsc{Theorem 12:} \textit{Consider equation (\ref{2nd kind eq}).\ Suppose
that Assumption 6 is satisfied, }$\Pr (g_{0}(c_{t})=0)=0$, \textit{for some }%
$w(Z_{t})$\textit{\ and }$\tilde{W}_{t}=(w(Z_{t}),c_{t})$\textit{,} $%
(c_{t+1},\tilde{W}_{t})$ \textit{is continuously distributed and there is
some }$j$ \textit{with }$A_{tj}=\delta _{0}R_{t+1,j}c_{t+1}^{-\gamma _{0}}$ 
\textit{satisfying }%
\begin{equation}
\text{E}[A_{tj}^{2}f_{c}(c_{t+1})^{-1}f_{\tilde{W}}(\tilde{W}_{t})^{-1}f_{c,%
\tilde{W}}(c_{t+1},\tilde{W}_{t})]<\infty .  \label{compac op}
\end{equation}%
\textit{Then (a) if }$\text{E}[A_{tj}\tilde{h}(c_{t+1},c_{t})|\tilde{W}%
_{t}]=0$ \textit{implies }$\tilde{h}(c_{t+1},c_{t})=0$\textit{\ a.s. and }$%
a(c_{t+1})+b(c_{t})=0\,$\textit{for }$c_{t}\in \mathcal{C}$ \textit{with }$%
\mathit{\Pr }(\mathcal{C)>}0$\textit{\ implies }$a(c_{t+1})$ \textit{is
constant then }$\Pi $ \textit{is nonsingular; (b) if }$g_{0}\in \mathcal{G}_{%
\bar{c}}\equiv \{g\in \mathcal{G}:g(\bar{c})\neq 0\}$ \textit{for some }$%
\bar{c}$\textit{\ and }$\text{E}[A_{tj}h(c_{t+1})|w(Z_{t}),c_{t}=\bar{c}]=0$ 
\textit{with }$h\in \mathcal{G}_{\bar{c}}$ \textit{implies }$h(c_{t+1})=0$%
\textit{\ a.s., then }$g_{0}$ \textit{is the unique solution to }$\text{E}%
[A_{t}g(c_{t+1})|W_{t}]=g(c_{t})$ \textit{up to scale.}

\bigskip

Equation (\ref{compac op}) implies $\text{E}[A_{tj}g(c_{t+1})|\tilde{W}_{t}]$
is a Hilbert-Schmidt integral operator and hence compact. Analogous
conditions could be imposed to ensure that $\mathcal{M}$ is closed. The
sufficient conditions for nonsingularity of $\Pi $ involve completeness of
the conditional expectation $\text{E}[A_{tj}h(c_{t+1},c_{t})|\tilde{W}_{t}]$
and a stronger version of a measurably separable condition from Florens,
Mouchart, and Rolin (1990). As previously noted, sufficient conditions for
completeness can be found in Newey and Powell (2003) and Andrews (2011) and
completeness is generic in the sense of Andrews (2011) and Lemma 3. A simple
sufficient condition for the measurably separable hypothesis is that the
support of $(c_{t+1},c_{t})$ is $\Re _{+}^{2}$, where $\Re _{+}=[0,\infty ).$

Condition b) is weaker than condition a). Condition (b) turns out to imply
global identification of $\beta _{0}$ and $g_{0}$ (up to scale) if $g_{0}(c)$
is bounded, and bounded away from zero. Because we focus on applying the
results of Section 4, we reserve this result to Theorem A.3 in Appendix D.
Even with global identification the result of Theorem 12 (a) is of interest,
because nonsingularity of $\Pi $ will be necessary for $\gamma _{0}$ to be
estimable at a root-n rate. The identification result for $\gamma _{0}$ in
Theorem A.3 involves large and small values of consumption growth, and so
amounts to identification at infinity, that may not lead to root-n
consistent estimation, e.g. see Chamberlain (1986).

A different approach to the nonparametric rank condition, that does not
require any instrument $w(Z_{t})$ in addition to $c_{t}$, can be based on
positivity of $g_{0}(c)$. The linear operator $\mathrm{E}%
[A_{tj}g(c_{t+1})|c_{t}]$ and $g(c)$ will be infinite dimensional
(functional) analogs of a positive matrix and a positive eigenvector
respectively, by equation (\ref{2nd kind eq}). The Perron-Frobenius Theorem
says that there is a unique positive eigenvalue and eigenvector (up to
scale) pair for a positive matrix. A functional analog, based on Krein and
Rutman (1950), gives uniqueness of $g_{0}(c),$ as well as of the discount
factor $\delta _{0}$. To describe this result let $r(c,s)=\mathrm{E}%
[R_{t+1,j}|c_{t+1}=s,c_{t}=c],$ $f(s,c)$ be the joint pdf of $%
(c_{t+1},c_{t}) $, $f(c)$ the marginal pdf of $c_{t}$ at $c$ and $%
K(c,s)=r(c,s)s^{-\gamma _{0}}f(s,c)/[f(s)f(c)].$ Then the equation $\mathrm{E%
}[A_{tj}g(c_{t+1})|c_{t}]=g(c_{t})$ can be written as%
\begin{equation}
\delta \int K(c,s)g(s)f(s)ds=g(c).  \label{eigen eq}
\end{equation}%
for $\delta =\delta _{0}\in (0,1]$. Here the matrix analogy is clear, with $%
K(c,s)f(s)$ being like a positive matrix, $g(c)$ an eigenvector, and $\delta
^{-1}$ an eigenvalue.

\bigskip

\textsc{Theorem 13:} \textit{Suppose that }$(R_{t,j},c_{t})$\textit{\ is
strictly stationary,\ $f(c,s)>0$ and $r(c,s)>0$ almost everywhere, and $\int
\int K(c,s)^{2}f(c)f(s)dcds<\infty .$ Then equation (\ref{eigen eq}) has a
unique positive solution $(\delta _{0},g_{0})$ in the sense that }$\delta
_{0}>0$\textit{, $g_{0}>0$ almost everywhere and }$E[g_{0}(c_{t})^{2}]=1$%
\textit{.}

\bigskip

The conditions of this result include $r(c,s)>0,$ that will hold if $%
R_{t+1,j}$ is a positive risk free rate. Under square-integrability of $K$,
we obtain global identification of the pair $(\delta _{0},g_{0})$. The
uniqueness of $g_{0}(c)$ in the conclusion of this result implies the
nonparametric rank condition. Note that by iterated expectation and
inclusion of $R_{t+1,j}$ in $R_{t+1}$ any solution to equation (\ref{2nd
kind eq}) must also satisfy equation (\ref{eigen eq}). Thus Theorem 13
implies that $g_{0}$ is the unique solution to (\ref{2nd kind eq}). Theorem
13 actually gives more, identification of the discount factor given
identification of $\gamma _{0}$.

Previously Escanciano and Hoderlein (2010) and Lewbel, Linton and Srisuma
(2012) considered nonparametric identification of marginal utility in
consumption level, $MU(C_{t})$, by solving the homogeneous integral equation
of the second kind:%
\begin{equation*}
\text{E}[R_{t+1,j}\delta _{0}MU(C_{t+1})|C_{t}]=MU(C_{t}).
\end{equation*}%
In particular, Escanciano and Hoderlein (2010) gave an insightful
identification result for the discount factor $\delta $ and the marginal
utility $MU$ based on the positivity of marginal utility and a version of
Perron-Frobenius theorem, but assuming that $C_{t}$ has a compact support,
that $MU$ is uniformly continuous on the support, and that $R_{t+1,j}$\ is a
risk-free rate. Lewbel, Linton and Srisuma (2012) used a genericity argument
for identification of $MU(C_{t})$. While we also use the positivity of the
unknown function $g$, we identify the \textquotedblleft correction term" $%
g(c)$ in the habit-formation model, rather than the marginal utility $MU$.
And we base our result on a Krein and Rutman (1950) theorem, a particular
functional analog of Perron-Frobenius, which allows us to completely avoid
making the compactness restriction on the support of $C_{t}$ and even of $%
c_{t}=C_{t+1}/C_{t}$.\footnote{%
The latter seems to be important to accommodate standard consumption-based
asset pricing models, in which $C_{t}$ is not stationary and has a
non-compact support of $(0,\infty )$.} Finally, Perron-Frobenius theory has
been extensively used by Hansen and Scheinkman (2009, 2012) and Hansen
(2012) in their research on the long-run risk and dynamic valuations in a
general Markov environment in which their valuation operators may not be
compact. Our results follow the same general path as all of the papers cited
in this paragraph, though we apply a Krein and Rutman (1950) theorem (and
its extensions) to a different operator than theirs.

The models considered here will generally be highly overidentified. We have
analyzed identification using only a single asset return $R_{t+1,j}.$ The
presence of more asset returns in equation (\ref{pri eqn}) provides
overidentifying restrictions. Also, in Theorem 12 we only use a function $%
w(Z_{t})$ of the available instrumental variables $Z_{t}$ in addition to $%
c_{t}$. The additional information in $Z_{t}$ may provide overidentifying
restrictions. These sources of overidentification are familiar in CCAPM
models. See, e.g., Hansen and Singleton (1982) and Chen and Ludvigson (2009).

\section{Conclusion}

We provide sufficient conditions for local identification for a general
class of semiparametric and nonparametric conditional moment restriction
models. We give new identification conditions for several important models
that illustrate the usefulness of our general results. In particular, we
provide primitive conditions for local identification in nonseparable
quantile IV models, single-index IV models, and semiparametric
consumption-based asset pricing models.

\appendix{}

\begin{center}
{\LARGE Appendix}
\end{center}

\section{Proofs for Section 2}

\setcounter{equation}{0} \renewcommand{\theequation}{A.\arabic{equation}}

\setcounter{theorem}{0} \renewcommand{\thetheorem}{A.\arabic{theorem}}

\subsection{Proof of Parametric Result}

By $rank(m^{\prime })=p,$ the nonnegative square root $\eta $ of the
smallest eigenvalue $\eta ^{2}$ of $(m^{\prime })^{T}m^{\prime }$ is
positive and $\left\vert m^{\prime }h\right\vert \geq \eta \left\vert
h\right\vert $ for $h\in \mathbb{R}^{p}.$ Also, by the definition of the
derivative there is $\varepsilon >0$ such that $\left\vert m(\alpha
)-m(\alpha _{0})-m^{\prime }(\alpha -\alpha _{0})\right\vert /\left\vert
\alpha -\alpha _{0}\right\vert <\eta $ for all $|\alpha -\alpha
_{0}|<\varepsilon $ with $\alpha \neq \alpha _{0}.$ Then%
\begin{equation}
\frac{\left\vert m(\alpha )-m^{\prime }(\alpha -\alpha _{0})\right\vert }{%
\left\vert m^{\prime }(\alpha -\alpha _{0})\right\vert }=\frac{\left\vert
m(\alpha )-m(\alpha _{0})-m^{\prime }(\alpha -\alpha _{0})\right\vert }{%
\left\vert \alpha -\alpha _{0}\right\vert }\frac{\left\vert \alpha -\alpha
_{0}\right\vert }{\left\vert m^{\prime }(\alpha -\alpha _{0})\right\vert }<%
\frac{\eta }{\eta }=1.  \label{ineq}
\end{equation}%
This inequality implies $m(\alpha )\neq 0,$ so $\alpha _{0}$ is locally
identified on $\{\alpha :|\alpha -\alpha _{0}|<\varepsilon \}$. $Q.E.D.$

\subsection{Proof of Theorem 1}

If $m^{\prime }h=m^{\prime }\tilde{h}$ for $h\neq \tilde{h}$ then for any $%
\lambda >0$ we have $m^{\prime }\bar{h}=0$ for $\bar{h}=\lambda (h-\tilde{h}%
)\neq 0$. For $\lambda $ small enough $\bar{h}$ would be in any open ball
around zero. Therefore, Assumption 1 holding on an open ball containing $%
\alpha _{0}$ implies that $m^{\prime }$ is invertible. By $m^{\prime }$ onto
and the Banach Inverse Theorem (Luenberger, 1969, p. 149) it follows that $%
(m^{\prime })^{-1}$ is continuous. Since any continuous linear map is
bounded, it follows that there exists $\eta >0$ such that $\left\Vert
m^{\prime }(\alpha -\alpha _{0})\right\Vert _{\mathcal{B}}\geq \eta
\left\Vert \alpha -\alpha _{0}\right\Vert _{\mathcal{A}}$ for all $\alpha
\in \mathcal{A}$.

Next, by Fr\'{e}chet differentiability at $\alpha _{0}$ there exists an open
ball $\mathit{\mathcal{N}}_{\varepsilon }$ centered at $\alpha _{0}$ such
that for all $\alpha \in \mathit{\mathcal{N}}_{\varepsilon }$, $\alpha \neq
\alpha _{0},$%
\begin{equation}
\frac{\left\Vert m(\alpha )-m(\alpha _{0})-m^{\prime }(\alpha -\alpha
_{0})\right\Vert _{\mathcal{B}}}{\left\Vert \alpha -\alpha _{0}\right\Vert _{%
\mathcal{A}}}<\eta .  \label{Frechet}
\end{equation}%
Therefore, at all such $\alpha \neq \alpha _{0}$,%
\begin{eqnarray*}
\frac{\left\Vert m(\alpha )-m^{\prime }(\alpha -\alpha _{0})\right\Vert _{%
\mathcal{B}}}{\left\Vert m^{\prime }(\alpha -\alpha _{0})\right\Vert _{%
\mathcal{B}}} &=&\frac{\left\Vert m(\alpha )-m(\alpha _{0})-m^{\prime
}(\alpha -\alpha _{0})\right\Vert _{\mathcal{B}}}{\left\Vert \alpha -\alpha
_{0}\right\Vert _{\mathcal{A}}}\frac{\left\Vert \alpha -\alpha
_{0}\right\Vert _{\mathcal{A}}}{\left\Vert m^{\prime }(\alpha -\alpha
_{0})\right\Vert _{\mathcal{B}}} \\
&<&\eta /\eta =1.
\end{eqnarray*}%
Therefore, as in the proof of the parametric result above, $m(\alpha )\neq 0$
for all $\alpha \in \mathit{\mathcal{N}}_{\varepsilon }$ with $\alpha \neq
\alpha _{0}$. \textit{Q.E.D.}

\subsection{Proof of Theorem 2}

Consider $\alpha \in $\textit{$\mathcal{N}$} with $\alpha \neq \alpha _{0}$.
Then 
\begin{equation*}
\frac{\left\Vert m(\alpha )-m^{\prime }(\alpha -\alpha _{0})\right\Vert _{%
\mathcal{B}}}{\left\Vert m^{\prime }(\alpha -\alpha _{0})\right\Vert _{%
\mathcal{B}}}=\frac{\left\Vert m(\alpha )-m(\alpha _{0})-m^{\prime }(\alpha
-\alpha _{0})\right\Vert _{\mathcal{B}}}{\left\Vert m^{\prime }(\alpha
-\alpha _{0})\right\Vert _{\mathcal{B}}}\leq \frac{L\left\Vert \alpha
-\alpha _{0}\right\Vert _{\mathcal{A}}^{r}}{\left\Vert m^{\prime }(\alpha
-\alpha _{0})\right\Vert _{\mathcal{B}}}<1.
\end{equation*}%
The conclusion follows as in the proof of Theorem 1. \textbf{\ }\textit{%
Q.E.D.}

\section{Proofs for Section 3}

\subsection{An Example for Lemma 3}

It is useful to give explicit examples of the randomization algorithms
obeying conditions 1 and 2 listed in Section 3. Suppose $\mathcal{A}=%
\mathcal{B}=L^{2}[0,1]$, and that $m^{\prime }$ is an integral operator 
\begin{equation*}
m^{\prime }\delta =\int K(\cdot ,t)\delta (t)dt.
\end{equation*}%
The kernel $K$ of this operator is generated as follows. The nature performs
step 1 by selecting two, possibly different, orthonormal bases $\{\phi
_{j}\} $ and $\{\varphi _{j}\}$ in $L^{2}[0,1]$. The nature performs step 2
by first selecting a bounded sequence $0<\sigma _{j}<\sigma $ for $j=0,1,...$%
, sampling $u_{j}$ as i.i.d. $U[-1,1]$, and then setting $\lambda
_{j}=u_{j}\sigma _{j}$. Finally, for some scalar $\kappa >0$ it sets 
\begin{equation*}
K=\kappa \left( \sum_{j=0}^{\infty }\lambda _{j}\phi _{j}\varphi _{j}\right)
.
\end{equation*}%
The operator defined in this way is well-defined over $\mathcal{A}$ and is
bounded, but it need not be compact. If compactness is required, we impose $%
\sum_{j=1}^{\infty }\sigma _{j}^{2}<\infty $ in the construction. If $K\geq
0 $ is required, we can impose $\phi _{0}=1$,$\varphi _{0}=1$, $|\varphi
_{j}|\leq c$ and $|\phi _{j}|\leq c$, for all $j$, where $c>1$ is a
constant, and $\sum_{j=0}^{\infty }\sigma _{j}<\infty $, and define instead $%
\lambda _{0}$ as $c\sum_{j=1}^{\infty }\lambda _{j}+|u_{0}|\sigma _{0}.$ If
in addition to positivity, $\int K(z,t)dt=1$ is required, for example if $%
K(z,t)=f(t|z)$ is a conditional density, then we select $\kappa >0$ so that $%
\kappa \lambda _{0}=1$. This algorithm for generating $m^{\prime }$
trivially obeys conditions 1 and 2 stated above. Furthermore, $u_{j}$ need
not be i.i.d. Take the extreme, opposite example, and set $u_{j}=U[-1,1]$
for all $j$, that is $u_{j}$'s are perfectly dependent. The resulting
algorithm for generating $m^{\prime }$ still trivially obeys conditions 1
and 2. The latter point -- of allowing perfect dependence -- is useful for
highlighting the differences with the approach and various examples given in
Andrews (2011); other than that, our point is the same.

An important example where dependence matters is the case with normal
instrumental regression, where the endogenous variable $X$ conditional on
the instrument $Z=z$ follows a normal distribution with mean $\rho z$ (and
variance normalized to 1). Here we let $\mathcal{A}=\mathcal{B}=L^{2}(%
\mathbb{R})$ equipped with standard normal density as a measure. In this
case, $m^{\prime }$ is an integral operator 
\begin{equation*}
m^{\prime }\delta =\int K(\cdot ,t)\delta (t)\frac{1}{\sqrt{2\pi }}%
e^{-t^{2}/2}dt,
\end{equation*}%
similarly to what we had above, where $K(t,z)$ has the following well known
representation: 
\begin{equation*}
K(t,z)=\sum_{j=0}^{\infty }\rho ^{j}\phi _{j}(t)\varphi _{j}(z),
\end{equation*}%
where $(\phi _{j})_{j=0}^{\infty }$ and $(\varphi _{j})_{j=0}^{\infty }$ are
the orthonormal (Hermite) polynomials. Hence if the nature draws $\rho $
from an absolutely continuous density on $(-1,1)$, then the full rank
condition holds with probability $1$. Note that the generalized Fourier
coefficients $(\rho ^{j})_{j=0}^{\infty }$ exhibit perfect dependence here.
To see that our randomization algorithm permits this, let the nature draw $%
\rho $ as specified above and draw $\lambda _{0}$ as an independent from $%
\rho $ random variable with support $(0,\infty )$ having an absolutely
continuous distribution. Then nature sets $\lambda _{j}=\rho ^{j}\lambda
_{0} $, $\kappa =1/\lambda _{0}$, and 
\begin{equation*}
K(t,z)=\kappa \left( \sum_{j=0}^{\infty }\lambda _{j}\phi _{j}(z)\varphi
_{j}(z)\right) .
\end{equation*}

\subsection{Proof of Lemma 3}

By assumptions there exists a compact, injective operator $\mathcal{K}:%
\mathcal{A}\mapsto \mathcal{B}$. By Theorem 15.16 in Kress (1999) $\mathcal{K%
}$ admits a singular value decomposition: 
\begin{equation*}
\mathcal{K}\delta =\sum_{j=0}^{N}\mu _{j}\langle \phi _{j},\delta \rangle
\varphi _{j},
\end{equation*}%
where $\{\phi _{j}\}$ is an orthonormal subset of $\mathcal{A}$, either
finite or countably infinite, with cardinality $N\leq \infty $, $\{\varphi
_{j}\}$ is an orthonormal subset of $\mathcal{B}$ of equal cardinality, and $%
\left( \mu _{j}\right) _{j=1}^{\infty }$ is bounded. Since $\Vert \mathcal{K}%
\delta \Vert _{\mathcal{B}}^{2}=\sum_{j=0}^{N}\mu _{j}^{2}\langle \phi
_{j},\delta \rangle ^{2},$ injectivity of $\mathcal{K}$ requires that $%
\{\phi _{j}\}$ must be an orthonormal basis in $\mathcal{A}$ and $\mu
_{j}\neq 0$ for all $j$. Therefore, step 1 is always feasible by using these 
$\{\phi _{j}\}$ and $\{\varphi _{j}\}$ in the construction. The order of
eigenvectors in these sets need not be preserved and could be arbitrary.
Step 2 is also feasible by using a product of Lebesgue- dominated measures
on a bounded subset of $\mathbb{R}$ to define a measure over $\mathbb{R}^{N}$%
, or, more generally, using any construction of measure on $\mathbb{R}^{N}$
from finite-dimensional measures obeying Kolmogorov's consistency conditions
(e.g. Dudley, 1989) and the additional condition that $\eta \{\lambda
_{j_{1}}\in A,\lambda _{j_{2}}\in \mathbb{R},...,\lambda _{j_{k}}\in \mathbb{%
R}\}=0$ if $\text{Leb}(A)=0$, for any finite subset $\{j_{1},...,j_{k}\}%
\subset \{0,...,N\}$. This verifies claim 1.

To verify claim 2, by Bessel's inequality we have that 
\begin{equation*}
\Vert m^{\prime }\delta \Vert _{\mathcal{B}}\geq \sum_{j=0}^{N}\lambda
_{j}^{2}\kappa ^{2}\langle \phi _{j},\delta \rangle ^{2}.
\end{equation*}%
$m^{\prime }$ is not injective iff $\lambda _{j}^{2}\kappa ^{2}=0$ for some $%
j$. By countable additivity and by $\text{Leb}(\{0\})=0\implies \eta
(\{\lambda _{j}=0\})=0$ holding by assumption, 
\begin{equation*}
\text{Pr}_{\eta }(\exists j\in \{0,1,...,N\}:\lambda _{j}=0)\leq
\sum_{j=0}^{N}\eta (\{\lambda _{j}=0\})=0.
\end{equation*}%
The final claim follows from the penultimate display. \emph{Q.E.D.} \newline

\subsection{Proof of Theorem 4}

By Assumption 3, for any $\alpha \neq \alpha _{0}$ and $\alpha \in \mathit{%
\mathcal{N}^{\prime \prime \prime }}$ with Fourier coefficients $b_{j}$ we
have%
\begin{equation*}
\left\Vert m^{\prime }(\alpha -\alpha _{0})\right\Vert _{\mathcal{B}}\geq
\left( \sum_{j}\mu _{j}^{2}b_{j}^{2}\right) ^{1/2}>L\left(
\sum_{j}b_{j}^{2}\right) ^{r/2}=L\left\Vert \alpha -\alpha _{0}\right\Vert _{%
\mathcal{A}}^{r},
\end{equation*}%
so the conclusion follows from Theorem 2. \textit{Q.E.D.}

\bigskip


\subsection{Proof of Corollary 5}

Consider $\alpha $ $\in \mathcal{N}^{\prime \prime \prime }$. Then 
\begin{equation}
\sum_{j}\mu _{j}^{-2/(r-1)}b_{j}^{2}<L^{-2/(r-1)}  \label{pf cor 4}
\end{equation}%
For $b_{j}=\left\langle \alpha -\alpha _{0},\phi _{j}\right\rangle $ note
that $\left\Vert \alpha -\alpha _{0}\right\Vert _{\mathcal{A}%
}=(\sum_{j}b_{j}^{2})^{1/2}$ by $\phi _{1},\phi _{2},...$ being an
orthonormal basis. Then 
\begin{eqnarray*}
(\sum_{j}b_{j}^{2})^{1/2} &=&\left( \sum_{j}\mu _{j}^{-2/r}\mu
_{j}^{2/r}b_{j}^{2}\right) ^{1/2}\leq \left( \sum_{j}\mu
_{j}^{-2/(r-1)}b_{j}^{2}\right) ^{(r-1)/2r}\left( \sum_{j}\mu
_{j}^{2}b_{j}^{2}\right) ^{1/2r} \\
&<&L^{-1/r}\left( \sum_{j}\mu _{j}^{2}b_{j}^{2}\right) ^{1/2r}\leq
L^{-1/r}\left( \left\Vert m^{\prime }(\alpha -\alpha _{0})\right\Vert _{%
\mathcal{B}}\right) ^{1/r},
\end{eqnarray*}%
where the first inequality holds by the H\"{o}lder inequality, the second by
eq. (\ref{pf cor 4}), and the third by Assumption 3. Raising both sides to
the $r^{th}$ power and multiplying through by $L$ gives%
\begin{equation}
L\left\Vert \alpha -\alpha _{0}\right\Vert _{\mathcal{A}}^{r}<\left\Vert
m^{\prime }(\alpha -\alpha _{0})\right\Vert _{\mathcal{B}}.  \label{key ineq}
\end{equation}%
The conclusion then follows from Theorem 4. \textit{Q.E.D.}

\bigskip

\subsection{Proof of Theorem 6}

Let $F(y|X,W)=\Pr (Y\leq y|X,W)$, $m(\alpha )=E\left[ 1(Y\leq \alpha (X))|W%
\right] -\tau $, and 
\begin{equation*}
m^{\prime }h=\text{E}\left[ f_{Y}(\alpha _{0}(X)|X,W)h(X)|W\right] ,
\end{equation*}%
so that by iterated expectations,%
\begin{equation*}
m(\alpha )=\text{E}\left[ F(\alpha (X)|X,W)|W\right] -\tau .
\end{equation*}%
\ Then by a pathwise mean value expansion, and by $f_{Y}(y|X,W)$
continuously differentiable%
\begin{eqnarray*}
&&\left\vert F(\alpha (X)|X,W)-F(\alpha _{0}(X)|X,W)-f_{Y}(\alpha
_{0}(X)|X,W)(\alpha (X)-\alpha _{0}(X))\right\vert \\
&=&\left\vert \left[ f_{Y}\left( \bar{\alpha}(X)|X,W\right) -f_{Y}(\alpha
_{0}(X)|X,W)\right] \left[ \alpha (X)-\alpha _{0}(X)\right] \right\vert \leq
L_{1}\left[ \alpha (X)-\alpha _{0}(X)\right] ^{2},
\end{eqnarray*}%
where $\bar{\alpha}(X)$ is the mean value of a pathwise Taylor expansion
that lies between $\alpha (X)$ and $\alpha _{0}(X)$. Then for $L_{1}L_{2}=L,$%
\begin{eqnarray*}
\left\vert m(\alpha )(W)-m(\alpha _{0})(W)-m^{\prime }(\alpha -\alpha
_{0})(W)\right\vert &\leq &L_{1}\text{E}\left[ \{\alpha (X)-\alpha
_{0}(X)\}^{2}|W\right] \\
&\leq &L\text{E}[\left\{ \alpha (X)-\alpha _{0}(X)\right\} ^{2}]=L\left\Vert
\alpha -\alpha _{0}\right\Vert _{\mathcal{A}}^{2}.
\end{eqnarray*}%
Therefore,%
\begin{equation*}
\left\Vert m(\alpha )-m(\alpha _{0})-m^{\prime }(\alpha -\alpha
_{0})\right\Vert _{\mathcal{B}}\leq L\left\Vert \alpha -\alpha
_{0}\right\Vert _{\mathcal{A}}^{2},
\end{equation*}%
so that Assumption 2 is satisfied with $r=2$ and $\mathcal{N}^{\prime \prime
}=\mathcal{A}$. The conclusion then follows from Corollary 5$.$ $Q.E.D.$

\bigskip

\section{Proofs for Section 4}

\subsection{Useful Results on Projections on Linear Subspaces}

Before proving the next Theorem we give two useful intermediate results. Let 
$\text{Proj}(b|\overline{\mathcal{M}})$ denote the orthogonal projection of
an element $b$ of a Hilbert space on a closed linear subset $\overline{%
\mathcal{M}}$ of that space.

\bigskip

\textsc{Lemma A.1: }\ \textit{If a) }$\overline{\mathcal{M}}$\textit{\ is a
closed linear subspace of a Hilbert space }$\mathcal{H};$ \textit{b) }$%
b_{j}\in \mathcal{H}\;(j=1,\ldots ,p)$\textit{; c) the }$p\times p$\textit{\
matrix }$\Pi $\textit{\ with }$\Pi _{jk}=\left\langle b_{j}-\text{Proj}%
(b_{j}|\overline{\mathcal{M}}),b_{k}-\text{Proj}(b_{k}|\overline{\mathcal{M}}%
)\right\rangle $\textit{\ is nonsingular, then for }$b=(b_{1},\ldots
,b_{p})^{T}$\textit{\ there exists }$\varepsilon >0$\textit{\ such that for
all }$a\in \mathbb{R}^{p}$\textit{\ and }$\zeta \in \overline{\mathcal{M}},$%
\begin{equation*}
\left\Vert b^{T}a+\zeta \right\Vert \geq \varepsilon \left( \left\vert
a\right\vert +\left\Vert \zeta \right\Vert \right) .
\end{equation*}

Proof: Let $\bar{b}_{j}=\text{Proj}(b_{j}|\overline{\mathcal{M}})$, $\tilde{b%
}_{j}=b_{j}-\bar{b}_{j},$ $\bar{b}=(\bar{b}_{1},...,\bar{b}_{p})^{T}$, and $%
\tilde{b}=(\tilde{b}_{1},...,\tilde{b}_{p})^{T}$. Note that for $\varepsilon
_{1}=\sqrt{\lambda _{\min }(\Pi )/2}$,%
\begin{eqnarray*}
\left\Vert b^{T}a+\zeta \right\Vert &=&\sqrt{\left\Vert \tilde{b}^{T}a+\zeta
+\bar{b}^{T}a\right\Vert ^{2}}=\sqrt{\left\Vert \tilde{b}^{T}a\right\Vert
^{2}+\left\Vert \zeta +\bar{b}^{T}a\right\Vert ^{2}} \\
&\geq &(\left\Vert \tilde{b}^{T}a\right\Vert +\left\Vert \zeta +\bar{b}%
^{T}a\right\Vert )/\sqrt{2}=(\sqrt{a^{T}\Pi a}+\left\Vert \zeta +\bar{b}%
^{T}a\right\Vert )/\sqrt{2} \\
&\geq &\varepsilon _{1}\left\vert a\right\vert +\left\Vert \zeta +\bar{b}%
^{T}a\right\Vert /\sqrt{2}.
\end{eqnarray*}%
Also note that for any $C^{\ast }\geq \sqrt{\sum_{j}\left\Vert \bar{b}%
_{j}\right\Vert ^{2}}$ it follows by the triangle and Cauchy-Schwartz
inequalities that%
\begin{equation*}
\left\Vert \bar{b}^{T}a\right\Vert \leq \sum_{j}\left\Vert \bar{b}%
_{j}\right\Vert \left\vert a_{j}\right\vert \leq C^{\ast }\left\vert
a\right\vert .
\end{equation*}%
Choose $C^{\ast }$ big enough that $\varepsilon _{1}/\sqrt{2}C^{\ast }\leq
1. $ Then by the triangle inequality,%
\begin{eqnarray*}
\left\Vert \zeta +\bar{b}^{T}a\right\Vert /\sqrt{2} &\geq &\left(
\varepsilon _{1}/\sqrt{2}C^{\ast }\right) \left\Vert \zeta +\bar{b}%
^{T}a\right\Vert /\sqrt{2}=\varepsilon _{1}\left\Vert \zeta +\bar{b}%
^{T}a\right\Vert /2C^{\ast } \\
&\geq &\varepsilon _{1}\left( \left\Vert \zeta \right\Vert -\left\Vert \bar{b%
}^{T}a\right\Vert \right) /2C^{\ast }\geq \varepsilon _{1}\left( \left\Vert
\zeta \right\Vert -C^{\ast }\left\vert a\right\vert \right) /2C^{\ast } \\
&=&\left( \varepsilon _{1}/2C^{\ast }\right) \left\Vert \zeta \right\Vert
-\varepsilon _{1}\left\vert a\right\vert /2.
\end{eqnarray*}%
Then combining the inequalities, for $\varepsilon =\min \{\varepsilon
_{1}/2,\varepsilon _{1}/2C^{\ast }\},$%
\begin{eqnarray*}
\left\Vert b^{T}a+\zeta \right\Vert &\geq &\varepsilon _{1}\left\vert
a\right\vert +\left( \varepsilon _{1}/2C^{\ast }\right) \left\Vert \zeta
\right\Vert -\varepsilon _{1}\left\vert a\right\vert /2 \\
&=&(\varepsilon _{1}/2)\left\vert a\right\vert +\left( \varepsilon
_{1}/2C^{\ast }\right) \left\Vert \zeta \right\Vert \geq \varepsilon \left(
\left\vert a\right\vert +\left\Vert \zeta \right\Vert \right) \text{. }Q.E.D.
\end{eqnarray*}

\bigskip

\textsc{Lemma A.2:} \textit{If Assumption 4 is satisfied then there is an }$%
\varepsilon >0$ \textit{such that for all }$(\beta ,g)\in \mathbb{R}%
^{p}\times \mathcal{N}_{g}^{\prime },$%
\begin{equation*}
\varepsilon (|\beta -\beta _{0}|+\left\Vert m_{g}^{\prime
}(g-g_{0})\right\Vert _{\mathcal{B}})\leq \left\Vert m^{\prime }(\alpha
-\alpha _{0})\right\Vert _{\mathcal{B}}.
\end{equation*}

Proof: Apply Lemma A.1 with $\mathcal{H}$ being the Hilbert space $\mathcal{B%
}$ described in Section 4, $\overline{\mathcal{M}}$ in Lemma A.1 being the
closed linear span of $\mathcal{M}\mathcal{=}\{m_{g}^{\prime }(g-g_{0}):g\in 
\mathcal{N}_{g}^{\prime }\},$ $b_{j}=m_{\beta }^{\prime }e_{j}$ for the $%
j^{th}$ unit vector $e_{j}$, and $a=\beta -\beta _{0}.$ Then for all $(\beta
,g)\in \mathbb{R}^{p}\times \mathcal{N}_{g}^{\prime }$ we have%
\begin{equation*}
m^{\prime }(\alpha -\alpha _{0})=b^{T}a+\zeta ,\text{ }b^{T}a=m_{\beta
}^{\prime }(\beta -\beta _{0}),\text{ }\zeta =m_{g}^{\prime }(g-g_{0})\in 
\overline{\mathcal{M}}\text{.}
\end{equation*}%
The conclusion then follows from the conclusion of Lemma A.1. \textit{Q.E.D.}

\bigskip

\subsection{Proof of Theorem 7}

Since Assumption 4 is satisfied the conclusion of Lemma A.2 holds. Let $%
\varepsilon $ be from the conclusion of Lemma A.2. Also let $\mathcal{N}_{g}=%
\mathcal{N}_{g}^{\prime }\cap \mathcal{N}_{g}^{\beta }$ for $\mathcal{N}%
_{g}^{\prime }$ from Assumption 4 and $\mathcal{N}_{g}^{\beta }$ from
Assumption 5. In addition let $B$ be from Assumption 5 with 
\begin{equation*}
\sup_{g\in \mathcal{N}_{g}^{\beta }}\text{E}[\sup_{\beta \in B}\left\vert
\partial \text{E}[\rho (Y,X,\beta ,g)|W]/\partial \beta -\partial \text{E}%
[\rho (Y,X,\beta _{0},g_{0})|W]/\partial \beta \right\vert ^{2}]<\varepsilon
^{2}.
\end{equation*}%
Then by $m(\beta _{0},g)$ linear in $g$ and expanding each element of $%
m(\beta ,g)(W)=\text{E}[\rho (Y,X,\beta ,g)|W]$ in $\beta ,$ it follows that
for each $(\beta ,g)\in B\times \mathcal{N}_{g}$, if $\beta \neq \beta _{0},$%
\begin{eqnarray*}
&&\left\Vert m(\alpha )-m^{\prime }(\alpha -\alpha _{0})\right\Vert _{%
\mathcal{B}}=\left\Vert m(\beta ,g)-m(\beta _{0},g)-m_{\beta }^{\prime
}(\beta -\beta _{0})\right\Vert _{\mathcal{B}} \\
&=&\left\Vert \left[ \partial m(\tilde{\beta},g)/\partial \beta -m_{\beta
}^{\prime }\right] (\beta -\beta _{0})\right\Vert _{\mathcal{B}}\leq
\left\Vert m_{\beta }^{\prime }(\tilde{\beta},g)-m_{\beta }^{\prime
}\right\Vert _{\mathcal{B}}\left\vert \beta -\beta _{0}\right\vert \\
&<&\varepsilon \left\vert \beta -\beta _{0}\right\vert \leq \varepsilon
(\left\vert \beta -\beta _{0}\right\vert +\left\Vert m_{g}^{\prime
}(g-g_{0})\right\Vert _{\mathcal{B}})\leq \left\Vert m^{\prime }(\alpha
-\alpha _{0})\right\Vert _{\mathcal{B}},
\end{eqnarray*}%
where $\tilde{\beta}$ is a mean value depending on $W$ that actually differs
from row to row of 
\begin{equation*}
m_{\beta }^{\prime }(\tilde{\beta},g)=\partial \text{E}[\rho (Y,X,\tilde{%
\beta},g)|W]/\partial \beta .
\end{equation*}%
Thus, $\left\Vert m(\alpha )-m^{\prime }(\alpha -\alpha _{0})\right\Vert _{%
\mathcal{B}}<\left\Vert m^{\prime }(\alpha -\alpha _{0})\right\Vert _{%
\mathcal{B}}$, implying $m(\alpha )\neq 0,$ giving the first conclusion.

To show the second conclusion, consider $(\beta ,g)\in \mathcal{N}$. If $%
\beta \neq \beta _{0}$ then it follows as above that $m(\alpha )\neq 0$. If $%
\beta =\beta _{0}$ and $g\neq g_{0}$ then by linearity in $g$ we have $%
\left\Vert m(\alpha )-m^{\prime }(\alpha -\alpha _{0})\right\Vert _{\mathcal{%
B}}=0$ while $\left\Vert m^{\prime }(\alpha -\alpha _{0})\right\Vert _{%
\mathcal{B}}=\left\Vert m_{g}^{\prime }(g-g_{0})\right\Vert _{\mathcal{B}%
}>0, $ so $m(\alpha )\neq 0$ follows as in the proof of Theorem 1. \textit{%
Q.E.D.}

\bigskip

\subsection{Proof of Theorem 8}

Since Assumption 4 is satisfied the conclusion of Lemma A.2 holds. Let $%
\varepsilon $ be from the conclusion of Lemma A.2. Define $B$ as in the
proof of Theorem 7. By Assumption 2, for $g\in \mathcal{N}_{g}^{\prime
\prime }$, $\left\Vert m(\beta _{0},g)-m_{g}^{\prime }(g-g_{0})\right\Vert _{%
\mathcal{B}}\leq L\left\Vert g-g_{0}\right\Vert _{\mathcal{A}}^{r}$. Then
similarly to the proof of Theorem 7 for all $\alpha \in \mathcal{N}$ with $%
\alpha \neq \alpha _{0}$,%
\begin{eqnarray*}
&&\left\Vert m(\alpha )-m^{\prime }(\alpha -\alpha _{0})\right\Vert _{%
\mathcal{B}} \\
&\leq &\left\Vert m(\beta ,g)-m(\beta _{0},g)-m_{\beta }^{\prime }(\beta
-\beta _{0})\right\Vert _{\mathcal{B}}+\left\Vert m(\beta
_{0},g)-m_{g}^{\prime }(g-g_{0})\right\Vert _{\mathcal{B}} \\
&<&\varepsilon \left\vert \beta -\beta _{0}\right\vert +L\left\Vert
g-g_{0}\right\Vert _{\mathcal{A}}^{r}\leq \varepsilon \left\vert \beta
-\beta _{0}\right\vert +\varepsilon \left\Vert m^{\prime
}(g-g_{0})\right\Vert _{\mathcal{B}} \\
&\leq &\left\Vert m^{\prime }(\alpha -\alpha _{0})\right\Vert _{\mathcal{B}}.
\end{eqnarray*}%
The conclusion follows as in the conclusion of Theorem 1. \textit{Q.E.D.}

\bigskip

\subsection{Proof of Theorem 9}

The proof will proceed by verifying the conditions of Theorem 7. Note that
Assumption 4 is satisfied. We now check Assumption 5. Note that for any $%
\delta >0$ and $g\in \mathcal{N}_{g}^{\delta }$, $g(X_{1}+X_{2}^{T}\beta )$
is continuously differentiable in $\beta $ with $\partial
g(X_{1}+X_{2}^{T}\beta )/\partial \beta =g^{\prime }(X_{1}+X_{2}^{T}\beta
)X_{2}.$ Also, for $\Delta $ a $p\times 1$ vector and $\bar{B}$ a
neighborhood of zero it follows by boundedness of $g_{0}^{\prime }$ and the
specification of $\mathcal{N}_{g}^{\delta }$ that for some $C>0$,%
\begin{equation*}
\text{E}[\sup_{\Delta \in \bar{B}}\left\vert g^{\prime
}(X_{1}+X_{2}^{T}(\beta +\Delta ))X_{2}\right\vert |W]\leq C\text{E}%
[\left\vert X_{2}\right\vert |W]<\infty \text{ a.s.}
\end{equation*}%
Therefore, by the dominated convergence theorem $m(\alpha )(W)=\text{E}%
[Y-g(X_{1}+X_{2}^{T}\beta )|W]$ is continuously differentiable in $\beta $
a.s. with 
\begin{equation*}
\partial m(\alpha )(W)/\partial \beta =-\text{E}[g^{\prime
}(X_{1}+X_{2}^{T}\beta )X_{2}|W].
\end{equation*}%
Next consider any $\varepsilon >0$ and let $B$ and $\delta $ satisfy%
\begin{equation*}
B=\{\beta :|\beta -\beta _{0}|^{2}<\varepsilon ^{2}/4C_{g}^{2}\text{E}%
[\left\vert X_{2}\right\vert ^{4}]\}\text{ and }\delta ^{2}<\varepsilon
^{2}/4\text{E}[\left\vert X_{2}\right\vert ^{2}].
\end{equation*}%
Then for $g\in \mathcal{N}_{g}^{\delta }$ we have, for $v(X,\beta
)=X_{1}+X_{2}^{T}\beta ,$%
\begin{eqnarray*}
&&\text{E}[\sup_{\beta \in B}\left\vert \partial m(\alpha )(W)/\partial
\beta -m_{\beta }^{\prime }(W)\right\vert ^{2}] \\
&=&\text{E}[\sup_{\beta \in B}\left\vert \text{E}[\{g^{\prime }(v(X,\beta
))-g_{0}^{\prime }(V)\}X_{2}|W]\right\vert ^{2}]\leq \text{E}[\left\vert
X_{2}\right\vert ^{2}\sup_{\beta \in B}\left\vert g^{\prime }(v(X,\beta
))-g_{0}^{\prime }(V)\right\vert ^{2}] \\
&\leq &2\text{E}[\left\vert X_{2}\right\vert ^{2}\sup_{\beta \in
B}\left\vert g^{\prime }(v(X,\beta ))-g_{0}^{\prime }(v(X,\beta
))\right\vert ^{2}]+2\text{E}[\left\vert X_{2}\right\vert ^{2}\sup_{\beta
\in B}\left\vert g_{0}^{\prime }(v(X,\beta ))-g_{0}^{\prime }(V)\right\vert
^{2}] \\
&\leq &2\delta ^{2}\text{E}[\left\vert X_{2}\right\vert ^{2}]+2C_{g}^{2}%
\text{E}[\left\vert X_{2}\right\vert ^{4}]\sup_{\beta \in B}\left\vert \beta
-\beta _{0}\right\vert ^{2}<\varepsilon ^{2}.
\end{eqnarray*}%
Thus Assumption 5 is satisfied so the first conclusion follows by the first
conclusion of Theorem 7. Also, $m_{g}^{\prime }(g-g_{0})=E[g(V)-g_{0}(V)|W]$
the rank condition for $m_{g}^{\prime }$ follows by the last bounded
completeness on $\mathcal{N}_{g}^{\prime }$, so that the final conclusion
follows by the final conclusion of Theorem 7. \textit{Q.E.D.}

\bigskip

\subsection{Proof of Theorem 10}

Suppose first that the conditional distribution of $W$ given $V$ is
complete. Note that by the projection definition, for all $h(V)$ with finite
mean-square we have%
\begin{equation*}
0=\text{E}[\{-\text{E}[g_{0}^{\prime }(V)X_{2j}|W]-\zeta _{j}^{\ast }(W)\}%
\text{E}[h(V)|W]]=\text{E}[\{-\text{E}[g_{0}^{\prime }(V)X_{2j}|W]-\zeta
_{j}^{\ast }(W)\}h(V)].
\end{equation*}%
Therefore, 
\begin{equation*}
\text{E}[-\text{E}[g_{0}^{\prime }(V)X_{2j}|W]-\zeta _{j}^{\ast }(W)|V]=0.
\end{equation*}%
Completeness of the conditional distribution of $W$ given $V$ then implies
that $-\text{E}[g_{0}^{\prime }(V)X_{2j}|W]-\zeta _{j}^{\ast }(W)=0,$ and
hence $\Pi _{jj}=0$. Since this is true for each $j$ we have $\Pi =0,$ $\Pi $
is singular.

Next, consider the second hypothesis and $\lambda \neq 0$. Let $\zeta
_{\lambda }^{\ast }(W)$ denote the projection of $-\text{E}[g_{0}^{\prime
}(V)\lambda ^{T}X_{2}|W]$ on $\overline{\mathcal{M}}$. Since $\text{E}%
[h(V)|W]=\text{E}[h(V)|T(W)]$ it follows that $\zeta _{\lambda }^{\ast }(W)$
is measurable with respect to (i.e. is a function of) $T(W).$ Since $\text{E}%
[g_{0}^{\prime }(V)\lambda ^{T}X_{2}|W]$ is not measurable with respect to $%
T(W)$, we have $-\text{E}[g_{0}^{\prime }(V)\lambda ^{T}X_{2}|W]-$ $\zeta
_{\lambda }^{\ast }(W)\neq 0,$ so that%
\begin{equation*}
\lambda ^{T}\Pi \lambda =\text{E}[\{-\text{E}[g_{0}^{\prime }(V)\lambda
^{T}X_{2}|W]-\zeta _{\lambda }^{\ast }(W)\}^{2}]>0.
\end{equation*}%
Since this is true for all $\lambda \neq 0$, it follows that $\Pi $ is
positive definite, and hence nonsingular. \textit{Q.E.D.}\bigskip

\section{Proofs for Section 5}

\subsection{Proof of Theorem 11}

The proof will proceed by verifying the conditions of Theorem 7 for $\rho
(Y,\beta ,g)$ from eq. (\ref{2nd kind eq}). We first check first part of
Assumption 4. Note that the mapping $m^{\prime }:\mathbb{R}^{2}\times 
\mathcal{G\longrightarrow B}$ is given by $m^{\prime }(\alpha -\alpha
_{0})=m_{\beta }^{\prime }(\beta -\beta _{0})+m_{g}^{\prime }(g-g_{0})$,
where%
\begin{eqnarray*}
m_{\beta }^{\prime }(\beta -\beta _{0}) &=&\text{E}%
[A_{t}g_{0}(c_{t+1})X_{t}^{T}|W_{t}](\beta -\beta _{0})\text{\quad and} \\
m_{g}^{\prime }(g-g_{0}) &=&\text{E}[A_{t}\{g(c_{t+1})-g_{0}(c_{t+1})%
\}|W_{t}]-\{g(c_{t})-g_{0}(c_{t})\}e.
\end{eqnarray*}%
Therefore the mapping $m^{\prime }$ is obviously linear. Since $\text{E}%
[D_{t}^{2}|W_{t}]$ and $\text{E}[D_{t}|W_{t}]$ exist with probability one by 
$\text{E}[D_{t}^{2}]<\infty $ and that $\left\vert A_{t}\right\vert ^{2}\leq
CD_{t}^{2}$. Then by the Cauchy-Schwartz inequality, for any $h\in \mathcal{G%
}$ we have by $D_{t}\geq 1\text{, E}[D_{t}^{2}|W_{t}]\geq 1$. 
\begin{eqnarray*}
&&\left\Vert \text{E}[A_{t}h(c_{t+1})|W_{t}]-h(c_{t})e\right\Vert _{\mathcal{%
B}}^{2}\leq C\text{E}[\text{E}[A_{t}^{T}h(c_{t+1})|W_{t}]\text{E}%
[A_{t}h(c_{t+1})|W_{t}]+h(c_{t})^{2}] \\
&\leq &C\text{E}[\text{E}[D_{t}^{2}|W_{t}]\text{E}[h(c_{t+1})^{2}|W_{t}]+C%
\text{E}[\text{E}[D_{t-1}^{2}|W_{t-1}]h(c_{t})^{2}]\leq C\left\Vert
h\right\Vert _{\mathcal{G}}^{2}.
\end{eqnarray*}%
Thus $m_{g}^{\prime }:\mathcal{G\longrightarrow B}$ is bounded. Also, noting
that $\left\vert m_{\beta }^{\prime }\right\vert \leq \text{E}%
[D_{t}g_{0}(c_{t+1})|W_{t}]$, the Cauchy-Schwartz inequality gives 
\begin{equation*}
\text{E}[\left\vert m_{\beta }^{\prime }(W)\right\vert ^{2}]\leq \text{E}[%
\text{E}[D_{t}^{2}|W_{t}]\text{E}[g_{0}(c_{t+1})^{2}|W_{t}]]\leq \left\Vert
g_{0}\right\Vert _{\mathcal{G}}^{2}<\infty
\end{equation*}%
and hence $m_{\beta }^{\prime }:\mathbb{R}^{2}\mathcal{\longrightarrow B}$
is bounded. Therefore the first part of Assumption 4 is satisfied with $%
\mathcal{N}_{g}^{\prime }=\mathcal{G}$.

Turning now to Assumption 5. Let $H_{t}(\beta ,g)=\delta
R_{t+1}c_{t+1}^{-\gamma }g(c_{t+1})$ and $B=[\delta _{0}-\Delta ,\delta
_{0}+\Delta ]\times \lbrack \gamma _{0}-\Delta ,\gamma _{0}+\Delta ]$. Note
that $H_{t}(\beta ,g)$ is twice continuously differentiable in $\beta $ and
by construction of $D_{t}$ that%
\begin{equation*}
\sup_{\beta \in B}\left\vert \frac{\partial H_{t}(\beta ,g)}{\partial \beta }%
\right\vert \leq D_{t}g(c_{t+1}),\text{ }\sup_{\beta \in B}\left\vert \frac{%
\partial ^{2}H_{t}(\beta ,g)}{\partial \beta _{j}\partial \beta }\right\vert
\leq D_{t}g(c_{t+1}),\text{ }(j=1,2).
\end{equation*}%
Therefore by standard results $\text{E}[\rho (Y_{t},\beta ,g)|W_{t}]=\text{E}%
[H_{t}(\beta ,g)|W_{t}]-g(c_{t})$ is twice continuously differentiable in $%
\beta $ on $B,$ $\partial \text{E}[\rho (Y_{t},\beta ,g)|W_{t}]/\partial
\beta =\text{E}[\partial H_{t}(\beta ,g)/\partial \beta |W_{t}].$ We also
have%
\begin{eqnarray*}
\left\vert \text{E}[\partial H_{t}(\beta ,g)/\partial \beta -\partial
H_{t}(\beta ,g_{0})/\partial \beta |W_{t}]\right\vert ^{2} &\leq &\text{E}%
[D_{t}^{2}|W_{t}]\text{E}[|g(c_{t+1})-g_{0}(c_{t+1})|^{2}|W_{t}], \\
\left\vert \text{E}[\partial H_{t}(\beta ,g_{0})/\partial \beta -\partial
H_{t}(\beta _{0},g_{0})/\partial \beta |W_{t}]\right\vert ^{2} &\leq &\text{E%
}[D_{t}^{2}|W_{t}]\text{E}[g_{0}(c_{t+1})^{2}|W_{t}]\left\vert \beta -\beta
_{0}\right\vert ^{2}.
\end{eqnarray*}%
Therefore we have%
\begin{eqnarray*}
&&\left\vert \frac{\partial \text{E}[\rho (Y,\beta ,g)|W]}{\partial \beta }-%
\frac{\partial \text{E}[\rho (Y,\beta _{0},g_{0})|W]}{\partial \beta }%
\right\vert ^{2} \\
&=&\left\vert \text{E}[\partial H_{t}(\beta ,g)/\partial \beta -\partial
H_{t}(\beta _{0},g_{0})/\partial \beta |W_{t}]\right\vert ^{2} \\
&\leq &2\text{E}[D_{t}^{2}|W_{t}]\{\text{E}%
[|g(c_{t+1})-g_{0}(c_{t+1})|^{2}|W_{t}]+\text{E}[g_{0}(c_{t+1})^{2}|W_{t}]%
\left\vert \beta -\beta _{0}\right\vert ^{2}\}.
\end{eqnarray*}%
Note that by iterated expectations, 
\begin{eqnarray*}
\text{E}[\text{E}[D_{t}^{2}|W_{t}]\text{E}%
[|g(c_{t+1})-g_{0}(c_{t+1})|^{2}|W_{t}]] &=&\left\Vert g-g_{0}\right\Vert _{%
\mathcal{G}}^{2}, \\
\text{E}[\text{E}[D_{t}^{2}|W_{t}]\text{E}[g_{0}(c_{t+1})|^{2}|W_{t}]]
&=&\left\Vert g_{0}\right\Vert _{\mathcal{G}}^{2}.
\end{eqnarray*}%
Consider any $\varepsilon >0$. Let 
\begin{equation*}
\mathcal{N}_{g}^{\beta }=\{g:\left\Vert g-g_{0}\right\Vert _{\mathcal{G}%
}\leq \varepsilon /2\}\text{ and }\tilde{B}=B\cap \{\beta :\left\vert \beta
-\beta _{0}\right\vert <\varepsilon /(2\left\Vert g_{0}\right\Vert _{%
\mathcal{G}}).
\end{equation*}%
Then for $g\in \mathcal{N}_{g}^{\beta }$ we have%
\begin{equation*}
\text{E}[\sup_{\beta \in \tilde{B}}\left\vert \partial m(\alpha
)(W)/\partial \beta -m_{\beta }^{\prime }(W)\right\vert ^{2}]\leq
2\left\Vert g-g_{0}\right\Vert _{\mathcal{G}}^{2}+2\left\Vert
g_{0}\right\Vert _{\mathcal{G}}^{2}\sup_{\beta \in \tilde{B}}\left\vert
\beta -\beta _{0}\right\vert ^{2}<\varepsilon ^{2}.
\end{equation*}%
Therefore Assumption 5 holds with $B$ there equal to $\tilde{B}$ here. The
conclusion then follows from Theorem 7. \textit{Q.E.D.}

\bigskip

\subsection{Proof of Theorem 12}

Let $\bar{a}(c_{t+1},\tilde{W}_{t})=\text{E}[A_{tj}|c_{t+1},\tilde{W}_{t}]$
and $\bar{d}(c_{t+1})=\text{E}[\text{E}[D_{t}^{2}|\tilde{W}_{t}]|c_{t+1}].$
Let $\mathcal{\tilde{B}=}\{b(\tilde{W}_{t}):\text{E}[b(\tilde{W}%
_{t})^{2}]<\infty \}$ and the operator $L:\mathcal{G\longrightarrow }%
\mathcal{\tilde{B}}$ be given by 
\begin{eqnarray*}
Lg &=&\text{E}[A_{tj}g(c_{t+1})|\tilde{W}_{t}]=\int \bar{a}(c,\tilde{W}%
_{t})g(c)\frac{f_{c,\tilde{W}}(c,\tilde{W}_{t})}{f_{\tilde{W}}(\tilde{W}_{t})%
}dc \\
&=&\int g(c)K(c,\tilde{W}_{t})f_{c}(c)\bar{d}(c)dc,\text{\quad }K(c,\tilde{W}%
_{t})=\frac{\bar{a}(c,\tilde{W}_{t})f_{c,\tilde{W}}(c,\tilde{W}_{t})}{f_{%
\tilde{W}}(\tilde{W}_{t})f_{c}(c)\bar{d}(c)}.
\end{eqnarray*}%
Note that $\bar{d}(c)\geq 1$ by $D_{t}^{2}\geq 1$. Therefore, 
\begin{eqnarray*}
\int K(c,w)^{2}\bar{d}(c)f_{c}(c)f_{\tilde{W}}(w)dcdw &=&\int \frac{\bar{a}%
(c,w)^{2}f_{c,\tilde{W}}(c,w)}{f_{\tilde{W}}(w)f_{c}(c)\bar{d}(c)}f_{c,%
\tilde{W}}(c,w)dcdw \\
&\leq &\int \frac{\bar{a}(c,w)^{2}f_{c,\tilde{W}}(c,w)}{f_{\tilde{W}%
}(w)f_{c}(c)}f_{c,\tilde{W}}(c,w)dcdw \\
&=&\text{E}[\text{E}[A_{tj}|c_{t+1},\tilde{W}_{t}]^{2}\frac{f_{c,\tilde{W}%
}(c_{t+1},\tilde{W}_{t})}{f_{c}(c_{t+1})f_{\tilde{W}}(\tilde{W}_{t})}] \\
&\leq &\text{E}[A_{tj}^{2}f_{c}(c_{t+1})^{-1}f_{\tilde{W}}(\tilde{W}%
_{t})^{-1}f_{c,\tilde{W}}(c_{t+1},\tilde{W}_{t})]<\infty .
\end{eqnarray*}%
It therefore follows by standard results that $L$ is Hilbert-Schmidt and
thus compact. Furthermore, it follows exactly as in the proof of Theorem 3.2
of Kress (1999), that 
\begin{equation*}
\mathcal{\tilde{M}=}\{\text{E}[A_{tj}g(c_{t+1})|\tilde{W}_{t}]-g(c_{t}):g\in 
\mathcal{G\}}
\end{equation*}%
is closed.

Next let $b=(b_{1},b_{2})^{T}$ be a constant vector and $\Delta
(c)=b_{1}/\delta _{0}-b_{2}\ln (c)$. Suppose $b^{T}\Pi b=0$. Then by the
definition of $\Pi $ there is $g_{k}\in \mathcal{G}$ such that 
\begin{equation*}
\text{E}[A_{t}g_{k}(c_{t+1})|\tilde{W}_{t}]-g_{k}(c_{t})e\longrightarrow 
\text{E}[A_{t}g_{0}(c_{t+1})\Delta (c_{t+1})|\tilde{W}_{t}]
\end{equation*}%
in mean square as $k\longrightarrow \infty $. It follows that for any $j$, 
\begin{equation*}
\text{E}[A_{tj}g_{k}(c_{t+1})|\tilde{W}_{t}]-g_{k}(c_{t})\longrightarrow 
\text{E}[A_{tj}g_{0}(c_{t+1})\Delta (c_{t+1})|\tilde{W}_{t}]
\end{equation*}%
in mean square. By $\mathcal{\tilde{M}}$ a closed set there exists $g^{\ast
}(c)$ such that%
\begin{equation}
\text{E}[A_{tj}g_{0}(c_{t+1})\Delta (c_{t+1})|\tilde{W}_{t}]=\text{E}%
[A_{tj}g^{\ast }(c_{t+1})|\tilde{W}_{t}]-g^{\ast }(c_{t})\text{.}
\label{gstar}
\end{equation}%
If $g^{\ast }(c_{t+1})=0$ then $\text{E}[A_{tj}g_{0}(c_{t+1})\Delta
(c_{t+1})|\tilde{W}_{t}]=0$ and by completeness of $\text{E}%
[A_{tj}h(c_{t+1},c_{t})|\tilde{W}_{t}]$ it follows that $g_{0}(c_{t+1})%
\Delta (c_{t+1})=0.$ Then by $\Pr (g_{0}(c_{t+1})\neq 0)=1$, we have $\Delta
(c_{t+1})=0$.

Next, suppose $\Pr (g^{\ast }(c_{t})\neq 0)>0.$ Then $\Pr (\min \{|g^{\ast
}(c_{t})|,g_{0}(c_{t})\}>0)>0,$ so for small enough $\varepsilon >0$ and $%
\mathcal{C=\{}c_{t}:\min \{|g^{\ast }(c_{t})|,g_{0}(c_{t})\}\geq \varepsilon
\}$ we have $\Pr (\mathcal{C})>0.$ Let $1_{t}^{\varepsilon }=$ $1(c_{t}\in 
\mathcal{C)}).$ Then multiplying through eq. (\ref{gstar}) by $%
1_{t}^{\varepsilon }/g^{\ast }(c_{t})$ and subtracting the conditional
expectation on the right-hand side gives 
\begin{equation*}
\text{E}[A_{tj}1_{t}^{\varepsilon }\frac{g_{0}(c_{t+1})\Delta
(c_{t+1})-g^{\ast }(c_{t+1})}{-g^{\ast }(c_{t})}|\tilde{W}%
_{t}]=1_{t}^{\varepsilon }\text{.}
\end{equation*}%
By eq. (\ref{pri eqn}) we also have%
\begin{equation*}
\text{E}[A_{tj}1_{t}^{\varepsilon }\left\{ \frac{g_{0}(c_{t+1})}{g_{0}(c_{t})%
}\right\} |\tilde{W}_{t}]=1_{t}^{\varepsilon }.
\end{equation*}%
By the completeness condition in part a) it then follows that%
\begin{equation*}
1_{t}^{\varepsilon }\frac{g_{0}(c_{t+1})\Delta (c_{t+1})-g^{\ast }(c_{t+1})}{%
-g^{\ast }(c_{t})}=1_{t}^{\varepsilon }\frac{g_{0}(c_{t+1})}{g_{0}(c_{t})}.
\end{equation*}%
Multiplying, dividing, and subtracting gives 
\begin{equation*}
1_{t}^{\varepsilon }\left[ \frac{g_{0}(c_{t+1})\Delta (c_{t+1})-g^{\ast
}(c_{t+1})}{-g_{0}(c_{t+1})}-\frac{g^{\ast }(c_{t})}{g_{0}(c_{t})}\right] =0.
\end{equation*}%
Then by the additive separability condition in part a) of the conditions it
follows that $g_{0}(c_{t+1})\Delta (c_{t+1})-g^{\ast
}(c_{t+1})=Cg_{0}(c_{t+1})$ for some nonzero constant $C$. Then by eq. (\ref%
{gstar}) and the second kind equation for $g_{0}$ we have%
\begin{equation*}
g^{\ast }(c_{t})=-\text{E}[A_{tj}\{g_{0}(c_{t+1})\Delta (c_{t+1})-g^{\ast
}(c_{t+1})\}|\tilde{W}_{t}]=-C\text{E}[A_{tj}g_{0}(c_{t+1})|\tilde{W}%
_{t}]=-Cg_{0}(c_{t})\text{.}
\end{equation*}%
Then by eq. (\ref{gstar}), 
\begin{equation*}
\text{E}[A_{tj}g_{0}(c_{t+1})\Delta (c_{t+1})|\tilde{W}_{t}]=\text{E}%
[A_{tj}g^{\ast }(c_{t+1})|\tilde{W}_{t}]-g^{\ast }(c_{t})=-C\{\text{E}%
[A_{tj}g_{0}(c_{t+1})|\tilde{W}_{t}]-g_{0}(c_{t})\}=0\text{.}
\end{equation*}%
By the completeness condition in part a) of the conditions it follows that $%
g_{0}(c_{t+1})\Delta (c_{t+1})=0$ so $\Delta (c_{t+1})=0$ follows by $\Pr
(g_{0}(c_{t})=0)=0$. Therefore, we find that $b^{T}\Pi b=0$ implies $\Delta
(c_{t+1})=0$. But we know that for $b\neq 0$ it is the case that $\Delta
(c_{t+1})\neq 0$. Therefore, $b\neq 0$ implies $b^{T}\Pi b>0$, i.e. $\Pi $
is nonsingular.

Next, under condition (b) of Theorem 12, if $\text{E}[A_{tj}g(c_{t+1})|%
\tilde{W}_{t}]=g(c_{t})$ for $g\in \mathcal{G}_{\bar{c}}$ it follows that
for $\bar{c}$ as given there,%
\begin{equation*}
\text{E}[A_{tj}\frac{g(c_{t+1})}{g(\bar{c})}|w(Z_{t}),c_{t}=\bar{c}]=1=\text{%
E}[A_{tj}\frac{g_{0}(c_{t+1})}{g_{0}(\bar{c})}|w(Z_{t}),c_{t}=\bar{c}].
\end{equation*}%
Then by the completeness condition in part (b) of the hypotheses, it follows
that $g(c_{t+1})/g(\bar{c})=g_{0}(c_{t+1})/g_{0}(\bar{c})$, i.e. 
\begin{equation*}
g(c_{t+1})=g_{0}(c_{t+1})g(\bar{c})/g_{0}(\bar{c}),
\end{equation*}%
so $g$ is equal to $g_{0}$ up to scale. This also implies that $g_{0}$ is
the unique solution to $\text{E}[A_{t}g(c_{t+1})|W_{t}]=g(c_{t})$ up to
scale. $Q.E.D.$

\subsection{Completeness and Global Identification in the CCAPM}

\textsc{T\textsc{heorem A.3.}} \textit{Consider model (\ref{pri eqn}). If }$%
(R_{t,j},c_{t})$ \textit{is strictly stationary}, $c_{t}$ \textit{is
continuously distributed with support }$[0,\infty )$\textit{,} $g_{0}(c)\geq
0$ \textit{is bounded and bounded away from zero, }$\mathrm{E}[\left\vert
R_{t+1,j}c_{t}^{-\gamma _{0}}\right\vert ]<\infty $, \textit{and there is }$%
\bar{c}$ \textit{such that }$\mathrm{E}[R_{t+1,j}h(c_{t+1})|w(Z_{t}),\bar{c}%
]=0$ \textit{and }$\mathrm{E}[\left\vert R_{t+1,j}h(c_{t+1})\right\vert
]<\infty $\textit{\ implies }$h(c_{t+1})=0$ \textit{then }$(\delta
_{0},\gamma _{0},g_{0})$ \textit{is identified (}$g_{0}$ \textit{up to
scale) among all }$(\delta ,\gamma ,g)$ \textit{with }$g(c)\geq 0,$ $g(c)$ 
\textit{bounded and bounded away from zero, and }$\mathrm{E}[\left\vert
R_{t+1,j}c_{t}^{-\gamma }\right\vert ]<\infty .$

\bigskip

Proof: Consider any two solutions $(\beta _{0},g_{0})$ and $(\beta
_{1},g_{1})$ to equation (\ref{pri eqn}) satisfying the conditions of
Theorem A.3. Then by iterated expectations, 
\begin{equation*}
\mathrm{E}\left[ R_{t+1,j}\delta _{0}c_{t+1}^{-\gamma _{0}}\frac{%
g_{0}(c_{t+1})}{g_{0}(\bar{c})}|w(Z_{t}),\bar{c}\right] =1=\mathrm{E}\left[
R_{t+1,j}\delta _{1}c_{t+1}^{-\gamma _{1}}\frac{g_{1}(c_{t+1})}{g_{1}(\bar{c}%
)}|w(Z_{t}),\bar{c}\right] .
\end{equation*}%
By completeness with $h(c_{t+1})=\delta _{0}c_{t+1}^{-\gamma
_{0}}g_{0}(c_{t+1})/g_{0}(\bar{c})-\delta _{1}c_{t+1}^{-\gamma
_{1}}g_{1}(c_{t+1})/g_{1}(\bar{c})$ it follows by multiplying and dividing
that 
\begin{equation*}
c_{t+1}^{\gamma _{1}-\gamma _{0}}=\frac{g_{1}(c_{t+1})}{g_{0}(c_{t+1})}\left[
\frac{\delta _{1}g_{0}(\bar{c})}{\delta _{0}g_{1}(\bar{c})}\right] .
\end{equation*}%
Since the object on the right is bounded and bounded away from zero and the
support of $c_{t+1}$ is $I=[0,\infty )$ it follows that $\gamma _{0}=\gamma
_{1}$. Then we have 
\begin{equation*}
g_{0}(c_{t+1})=g_{1}(c_{t+1})\left[ \frac{\delta _{1}g_{0}(\bar{c})}{\delta
_{0}g_{1}(\bar{c})}\right] \text{ \ }a.e.\text{ in }I^{2},
\end{equation*}%
so that there is a constant $D>0$ such that $g_{0}(c_{t+1})=Dg_{1}(c_{t+1})$
a.e. in $I$. We can also assume that $g_{0}(\bar{c})=Dg_{1}(\bar{c})$ since $%
c_{t}$ is continuously distributed. Substituting then gives $D=(\delta
_{1}/\delta _{0})D$, implying $\delta _{1}=\delta _{0}$. \textit{Q.E.D.}

\bigskip

Previously Chen and Ludvigson (2009) show global identification of $(\delta
_{0},\gamma _{0},g_{0})$\ under different conditions. In their results $%
\mathrm{E}\left[ R_{t+1,j}h(c_{t+1},c_{t})|w(Z_{t}),c_{t}\right] $ is
assumed to be complete, which is similar to condition (a) in Theorem 12 and
is stronger than completeness at $c_{t}=\bar{c}$, but $g(c)$ is not assumed
to be bounded or bounded away from zero on $[0,\infty )$.

\subsection{A Useful Result on Uniqueness and Existence of Positive
Eigenfunctions}

The following result and its proof in part rely on the fundamental results
of Krein and Rutman (1950), specifically their Theorem 6.1 and example $%
\beta ^{\prime }$. Krein and Rutman (1950) is one of many extensions of the
Perron-Frobenius theory of positive matrices to the case of operators
leaving invariant a cone in a Banach space.

Let $I$ be a Borel subset of $\mathbb{R}^m$ and $\mu$ be a $\sigma$-finite
measure with support $I$. Consider the space $L^2(\mu)$, equipped with the
standard norm $\| \cdot\|$. We consider the following conditions on the
kernel $K$:

\begin{itemize}
\item[1.] $K(s,t)$ is a non-negative, measurable kernel such that $\int \int
K^2(s,t) d\mu(t) d \mu(s) < \infty$.

\item[2.] $K(s,t) = 0$ on a set of points $(t,s)$ of measure zero under $\mu
\times \mu$.
\end{itemize}

Consider an integral operator $L$ from $L^{2}(\mu )$ to $L^{2}(\mu )$
defined by: 
\begin{equation*}
L\varphi :=\int K(\cdot ,t)\varphi (t)d\mu (t),
\end{equation*}%
and its adjoint operator 
\begin{equation*}
L^{\ast }\psi :=\int K(t,\cdot )\psi (t)d\mu (t).
\end{equation*}%
It is known that these operators are compact under condition 1. The lemma
given below shows that under these assumptions we have existence and global
uniqueness of the positive eigenpair $(\rho ,\varphi )$ such that $L\varphi
=\rho \varphi $, in the sense that is stated below. This lemma extends
example $\beta ^{\prime }$ outlined in Krein and Rutman (1950) that looked
at the complex Hilbert space $L^{2}[a,b]$, $0<a<b<\infty $, an extension
which we were not able to track easily in the literature, so we simply
derived it; we also provided an additional step (3), not given in the
outline, to fully verify uniqueness. Note that we removed the complex
analysis based arguments, since they are not needed here.

\bigskip

\textsc{Lemma A.4.} \textit{Under conditions 1 and 2, there exists a unique
eigenpair $(\rho ,\varphi )$, consisting of an eigenvalue $\rho $ and
eigenfunction $\varphi $ such that $L\varphi =\rho \varphi $ and $\rho >0$, $%
\Vert \varphi \Vert =1$, $\varphi \geq 0$; moreover, $\varphi >0$ $\mu $-a.e.%
}

\bigskip

Proof. The proof is divided in five steps.

\medskip

(1) Let $C^{o}$ be the cone of nonnegative functions in $A=L^{2}(\mu )$. In
the proof we shall use the following result on the existence of non-negative
eigenpair from Krein and Rutman (1950, Theorem 6.1).

\begin{quote}
Consider a cone $C^o$ in a Banach space $A$ such that the closure of the
linear hull of $C^o$ is $A$. Consider a linear, compact operator $L: A
\mapsto A$ such that $L C^o \subset C^o$, and that has one point of spectrum
different from zero. Then it has a positive eigenvalue $\rho$, not less in
modulus than every other eigenvalues, and to this eigenvalue there
corresponds at least one eigenvector $\varphi \in C^o$ of the operator $L$ ($%
L \varphi = \rho \varphi$) and at least one eigenvector $\psi \neq 0$ of the
dual operator $L^*$ ($L^* \psi = \rho \psi$).
\end{quote}

The theorem requires that the closure of the linear hull of the cone is $A$.
This is true in our case for $A = L^2(\mu)$ and the cone $C^o$ of the
non-negative functions in $A$, since $C^o - C^o$ is dense in $A$. Moreover,
since 
\begin{equation*}
\sigma_2 = \int K(s,t) K(t,s) d \mu(s) d \mu(t)>0,
\end{equation*}
which is equal to sum of squared eigenvalues of $L$, the spectrum of $L$
must have at least one point different from zero. Therefore, application of
the theorem quoted above implies that there exists $\rho>0$ and $\varphi$
and $\psi $ s.t. $\mu$-a.e. 
\begin{eqnarray}
& & \varphi (s) = \rho^{-1} \int K(s,t) \varphi(t) d \mu(t), \ \ \varphi
\geq 0, \| \varphi\| = 1, \ \rho>0;  \label{eq: primary equation} \\
& & \psi(s) = \rho^{-1} \int K(t,s) \psi(t) d \mu(t), \ \ \| \psi \| =1.
\label{eq: dual equation}
\end{eqnarray}

(2) We would like to prove that any eigenvalue $\rho>0 $ associated to a
nonnegative eigenfunction $\varphi \geq 0$ must be a simple eigenvalue, i.e. 
$\varphi$ is the only eigenfunction in $L^2 (\mu)$ associated with $\rho$.
For this purpose we shall use the following standard fact on linear compact
operators, e.g. stated in Krein and Rutman (1950) and specialized to our
context: An eigenvalue $\rho$ of $L$ is simple if and only if the equations $%
L \varphi = \rho \varphi$ and $L^* \psi = \rho \psi$ have no solutions
orthogonal to each other, i.e. satisfying $\varphi \neq 0$, $\psi \neq 0$, $%
\int \psi(s) \varphi(s) d \mu(s) = 0$. So for this purpose we will show in
steps (4) and (5) below that $\psi$ is of constant sign $\mu$-a.e. and $%
\varphi$ and $\psi$ only vanish on a set of measure $0$ under $\mu$. Since $%
\varphi \geq 0$, this implies 
\begin{equation*}
\int \psi(s) \varphi(s) d \mu(s) \neq 0,
\end{equation*}
and we conclude from the quoted fact that $\rho$ is a simple eigenvalue. 
\newline

(3) To assert the uniqueness of the nonnegative eigenpair $(\rho, \varphi)$
(meaning that $L \varphi = \rho \varphi$, $\rho>0$, $\varphi \geq 0$, $%
\|\varphi\|=1$), suppose to the contrary that there is another nonnegative
eigenpair $(r, \zeta)$ . Then $r$ is also an eigenvalue of $L^*$ by the
Fredholm theorem (Kress, 1999, Theorem 4.14), which implies by definition of
the eigenvalue that there exists a dual eigenfunction $\eta \neq 0$ such
that $L^* \eta = r\eta$ and $\|\eta \|=1$.

By step (4) below we must have $\zeta>0$, $\varphi > 0$ $\mu$-a.e. Hence by
step (5) the dual eigenfunctions $\eta$ and $\psi$ are non-vanishing and of
constant sign $\mu$-a.e., which implies $\int \eta(s) \varphi(s) d \mu(s)
\neq 0$. Therefore, $r= \rho$ follows from the equality: 
\begin{eqnarray*}
r \int \eta(s) \varphi(s) d \mu(s) = \int \int K(t,s) \eta(t) d \mu(t)
\varphi(s) d \mu(s) = \rho \int \eta(t) \varphi(t) d \mu(t).
\end{eqnarray*}

(4) Let us prove that any eigenfunction $\varphi \geq 0$ of $L$ associated
with an eigenvalue $\rho>0$ must be $\mu$-a.e. positive. Let $S$ denote the
set of zeros of $\varphi$. Evidently, $\mu(S) < \mu(I)$. If $s \in S$, then 
\begin{equation*}
\int K(s,t) \varphi(t) d \mu(t) =0.
\end{equation*}
Therefore $K(s,t)$ vanishes almost everywhere on $(s,t) \in S \times
(I\setminus S)$. However the set of zeroes of $K(s,t)$ is of measure zero
under $\mu \times \mu$, so $\mu ( S) \times \mu (I\setminus S)=0$, implying $%
\mu (S) =0$. \newline

(5) Here we show that any eigen-triple $(\rho ,\varphi ,\psi )$ solving (\ref%
{eq: dual equation}) and (\ref{eq: primary equation}) obeys: 
\begin{equation}
\text{sign}(\psi (s))=1\text{ $\mu $-a.e. or }\text{sign}(\psi (s))=-1\text{ 
$\mu $-a.e.}  \label{eq: sign equality}
\end{equation}

From equation (\ref{eq: dual equation}) it follows that $\mu $-a.e. 
\begin{equation*}
|\psi (s)|\leq \rho ^{-1}\int K(t,s)|\psi (t)|d\mu (t).
\end{equation*}%
Multiplying both sides by $\varphi (s)$, integrating and applying (\ref{eq:
primary equation}) yields 
\begin{equation*}
\int |\psi (s)|\varphi (s)d\mu (s)\leq \rho ^{-1}\int \int K(t,s)\varphi
(s)|\psi (t)|d\mu (t)d\mu (s)=\int |\psi (t)|\varphi (t)d\mu (t).
\end{equation*}%
It follows that $\mu $-a.e. 
\begin{equation*}
|\psi (s)|=\rho ^{-1}\int K(t,s)|\psi (t)|d\mu (t),
\end{equation*}%
i.e. $|\psi |$ is an eigenfunction of $L^{\ast }$.

Next, equation $|\psi(s)| = \psi(s) \text{sign}(\psi(s) )$ implies that $\mu$%
-a.e. 
\begin{equation*}
\rho^{-1} \int K(t,s) | \psi(t) | d \mu(t) = \rho^{-1} \int K(t,s) \psi(t) d
\mu(t) \text{sign}(\psi(s)).
\end{equation*}
It follows that for a.e. $(t,s)$ under $\mu \times \mu$ 
\begin{equation*}
|\psi(t)| = \psi(t) \text{sign}(\psi(s)).
\end{equation*}
By the positivity condition on $K$, $|\psi|>0$ $\mu$-a.e. by the same
reasoning as given in step (4). Thus, (\ref{eq: sign equality}) follows. 
\textit{Q.E.D.}\newline

\subsection{Proof of Theorem 13.}

Note that $K(c,s)=r(c,s)s^{-\gamma _{0}}f(s,c)/[f(s)f(c)]>0$ almost
everywhere by $r(c,s)>0$ and $f(s,c)>0$ almost everywhere. Therefore the
conclusion follows from Lemma A.4 with $f(s)ds=d\mu (s)$. \textit{Q.E.D.}

\section{Tangential Cone Conditions}

In this Appendix we discuss some inequalities that are related to
identification of $\alpha _{0}$. Throughout this Appendix we maintain that $%
m(\alpha _{0})=0$. Define%
\begin{eqnarray*}
\mathcal{N} &=&\{\alpha :m(\alpha )\neq 0\},\quad \mathcal{N}^{\prime
}=\{\alpha :m^{\prime }(\alpha -\alpha _{0})\neq 0\}, \\
\mathcal{N}_{\eta }^{\prime } &=&\{\alpha :\left\Vert m(\alpha )-m^{\prime
}(\alpha -\alpha _{0})\right\Vert _{\mathcal{B}}\leq \eta \left\Vert
m^{\prime }(\alpha -\alpha _{0})\right\Vert _{\mathcal{B}}\},\eta >0, \\
\mathcal{N}_{\eta } &=&\{\alpha :\left\Vert m\left( \alpha \right)
-m^{\prime }\left( \alpha -\alpha _{0}\right) \right\Vert _{\mathcal{B}}\leq
\eta \left\Vert m\left( \alpha \right) \right\Vert _{\mathcal{B}}\},\quad
\eta >0.
\end{eqnarray*}%
Here $\mathcal{N}$ can be interpreted as the identified set and $\mathcal{N}%
^{\prime }$ as the set where the rank condition holds. The set $\mathcal{N}%
_{\eta }^{\prime }$ is a set on which an inequality version of equation (\ref%
{tang cone cond}) holds. The inequality used to define $\mathcal{N}_{\eta }$
is similar to the tangential cone condition from the literature on
computation in nonlinear ill-posed inverse problems, e.g. Hanke, Neubauer
and Scherzer (1995) and Dunker et. al. (2012).

The following results gives some relations among these sets:

\bigskip

\textsc{Lemma A.5: }\textit{For any }$\eta >0$\textit{,} 
\begin{equation*}
\mathcal{N}_{\eta }\cap \mathcal{N}^{\prime }\subset \mathcal{N},\quad 
\mathcal{N}_{\eta }^{\prime }\cap \mathcal{N}\subset \mathcal{N}^{\prime }.
\end{equation*}%
\textit{If }$0<\eta <1$\textit{\ then}%
\begin{equation*}
\mathcal{N}_{\eta }\cap \mathcal{N}\subset \mathcal{N}^{\prime },\quad 
\mathcal{N}_{\eta }^{\prime }\cap \mathcal{N}^{\prime }\subset \mathcal{N}.
\end{equation*}%
Proof: Note that $\alpha \in \mathcal{N}_{\eta }$ and the triangle
inequality gives%
\begin{equation*}
-\left\Vert m(\alpha )\right\Vert _{\mathcal{B}}+\left\Vert m^{\prime
}(\alpha -\alpha _{0})\right\Vert _{\mathcal{B}}\leq \eta \left\Vert
m(a)\right\Vert _{\mathcal{B}}
\end{equation*}%
so that $\left\Vert m(\alpha )\right\Vert _{\mathcal{B}}\geq (1+\eta
)^{-1}\left\Vert m^{\prime }(\alpha -\alpha _{0})\right\Vert _{\mathcal{B}}.$
Therefore if $\alpha \in \mathcal{N}_{\eta }\cap \mathcal{N}^{\prime }$ we
have $\left\Vert m(\alpha )\right\Vert _{\mathcal{B}}>0$, i.e. $\alpha \in 
\mathcal{N}$, giving the first conclusion. \ Also, if $\alpha \in \mathcal{N}%
_{\eta }^{\prime }$ we have 
\begin{equation*}
-\left\Vert m^{\prime }(\alpha -\alpha _{0})\right\Vert _{\mathcal{B}%
}+\left\Vert m(\alpha )\right\Vert _{\mathcal{B}}\leq \eta \left\Vert
m^{\prime }(\alpha -\alpha _{0})\right\Vert _{\mathcal{B}}
\end{equation*}%
so that $\left\Vert m^{\prime }\left( \alpha -\alpha _{0}\right) \right\Vert
_{\mathcal{B}}\geq (1+\eta )^{-1}\left\Vert m(\alpha )\right\Vert _{\mathcal{%
B}}.$ Therefore, if $\alpha \in \mathcal{N}_{\eta }^{\prime }\cap \mathcal{N}
$ we have $\left\Vert m^{\prime }(\alpha -\alpha _{0})\right\Vert _{\mathcal{%
B}}>0$, giving the second conclusion.

Next, for $0<\eta <1$ and $\alpha \in \mathcal{N}_{\eta }$ we have%
\begin{equation*}
\left\Vert m(\alpha )\right\Vert _{\mathcal{B}}-\left\Vert m^{\prime
}(\alpha -\alpha _{0})\right\Vert _{\mathcal{B}}\leq \eta \left\Vert
m(\alpha )\right\Vert _{\mathcal{B}}
\end{equation*}%
so that $\left\Vert m^{\prime }(\alpha -\alpha _{0})\right\Vert \geq (1-\eta
)\left\Vert m(\alpha )\right\Vert _{\mathcal{B}}.$ \ Therefore, if $\alpha
\in \mathcal{N}_{\eta }\cap \mathcal{N}$ we have $\left\Vert m^{\prime
}(\alpha -\alpha _{0})\right\Vert _{\mathcal{B}}>0$, giving the third
conclusion. \ Similarly, for $0<\eta <1$ and $\alpha \in \mathcal{N}_{\eta
}^{\prime }$ we have%
\begin{equation*}
\left\Vert m^{\prime }(\alpha -\alpha _{0})\right\Vert _{\mathcal{B}%
}-\left\Vert m(\alpha )\right\Vert _{\mathcal{B}}\leq \eta \left\Vert
m^{\prime }(\alpha -\alpha _{0})\right\Vert _{\mathcal{B}}
\end{equation*}%
so that $\left\Vert m(\alpha )\right\Vert _{\mathcal{B}}\geq (1-\eta
)\left\Vert m^{\prime }(\alpha -\alpha _{0})\right\Vert _{\mathcal{B}}$.
Therefore if $\alpha \in \mathcal{N}_{\eta }^{\prime }\cap \mathcal{N}%
^{\prime }$ we have $\left\Vert m(\alpha )\right\Vert _{\mathcal{B}}>0$,
giving the fourth conclusion. \textit{Q.E.D.}

\bigskip

The first conclusion shows that when the tangential cone condition is
satisfied the set on which the rank condition holds is a subset of the
identified set. The second condition is less interesting, but does show that
the rank condition is necessary for identification when $\alpha \in \mathcal{%
N}_{\eta }^{\prime }$. The third conclusion shows that the rank condition is
also necessary for identification under the tangential cone condition for $%
0<\eta <1.$ The last conclusion shows that when $\alpha \in \mathcal{N}%
_{\eta }^{\prime }$ with $0<\eta <1$ the rank condition is sufficient for
identification.

When the side condition that $\alpha \in \mathcal{N}_{\eta }$ or $\alpha \in 
\mathcal{N}_{\eta }^{\prime }$ are imposed for $0<\eta <1$, the rank
condition is necessary and sufficient for identification.

\bigskip

\textsc{Corollary A.6}: \textit{If }$0<\eta <1$\textit{\ then}%
\begin{equation*}
\mathcal{N}_{\eta }\cap \mathcal{N}^{\prime }=\mathcal{N}_{\eta }\cap 
\mathcal{N},\quad \mathcal{N}_{\eta }^{\prime }\cap \mathcal{N}^{\prime }=%
\mathcal{N}_{\eta }^{\prime }\cap \mathcal{N}.
\end{equation*}

\bigskip

Proof: By intersecting both sides of the first conclusion of Lemma A.5 with $%
\mathcal{N}_{\eta }$ we find that $\mathcal{N}_{\eta }\cap \mathcal{N}%
^{\prime }\subset \mathcal{N}_{\eta }\cap \mathcal{N}$. \ For $\eta <1$ it
follows similarly from the third conclusion of Lemma A.5 that $\mathcal{N}%
_{\eta }\cap \mathcal{N}\subset \mathcal{N}_{\eta }\cap \mathcal{N}^{\prime
} $, implying $\mathcal{N}_{\eta }\cap \mathcal{N}^{\prime }=\mathcal{N}%
_{\eta }\cap \mathcal{N},$ the first conclusion. \ The second conclusion
follows similarly. \textit{Q.E.D.}

\bigskip

The equalities in the conclusion of this result show that the rank condition
(i.e. $\alpha \in \mathcal{N}^{\prime }$) is necessary and sufficient for
identification (i.e. $\alpha \in \mathcal{N}$) under either of the side
conditions that 
\begin{equation*}
\alpha \in \mathcal{N}_{\eta }^{\prime }\text{ or }\alpha \in \mathcal{N}%
_{\eta },\text{ }0<\eta <1.
\end{equation*}%
In parametric models Rothenberg (1971) showed that when the Jacobian has
constant rank in a neighborhood of the true parameter the rank condition is
necessary and sufficient for local identification. These conditions fill an
analogous role here, in the sense that when $\alpha $ is restricted to
either set, the rank condition is necessary and sufficient for
identification.

The sets $\mathcal{N}_{\eta }$ and $\mathcal{N}_{\eta }^{\prime }$ are
related to each other in the way shown in the following result.

\bigskip

\textsc{Lemma A.7.} \textit{If }$0<\eta <1$ \textit{\ then } $\mathcal{N}
_{\eta }\subset $\textit{\ }$\mathcal{N}_{\eta /(1-\eta )}^{\prime }$\textit{%
\ and }$\mathcal{N}_{\eta }^{\prime }\subset $\textit{\ }$\mathcal{N}_{\eta
/(1-\eta )}$\textit{. }

\bigskip

Proof: By the triangle inequality%
\begin{eqnarray*}
\left\Vert m^{\prime }(\alpha -\alpha _{0})\right\Vert _{\mathcal{B}} &\leq
&\left\Vert m(\alpha )-m^{\prime }(\alpha -\alpha _{0})\right\Vert _{%
\mathcal{B}}+\left\Vert m(\alpha )\right\Vert _{\mathcal{B}}, \\
\left\Vert m(\alpha )\right\Vert _{\mathcal{B}} &\leq &\left\Vert m(\alpha
)-m^{\prime }(\alpha -\alpha _{0})\right\Vert _{\mathcal{B}}+\left\Vert
m^{\prime }(\alpha -\alpha _{0})\right\Vert _{\mathcal{B}}.
\end{eqnarray*}%
Therefore, for $\alpha \in \mathcal{N}_{\eta },$ 
\begin{equation*}
\left\Vert m(\alpha )-m^{\prime }(\alpha -\alpha _{0})\right\Vert _{\mathcal{%
B}}\leq \eta \left\Vert m(\alpha )-m^{\prime }(\alpha -\alpha
_{0})\right\Vert _{\mathcal{B}}+\eta \left\Vert m^{\prime }(\alpha -\alpha
_{0})\right\Vert _{\mathcal{B}}.
\end{equation*}%
Subtracting $\eta \left\Vert m(\alpha )-m^{\prime }(\alpha -\alpha
_{0})\right\Vert _{\mathcal{B}}$ from both sides and dividing by $1-\eta $
gives $\alpha \in $ $\mathcal{N}_{\eta /(1-\eta )}^{\prime }$. The second
conclusion follows similarly. \textit{Q.E.D.}

\end{document}